\documentclass[12pt,verbatim]{amsart}
\usepackage{tabularx}
\usepackage{longtable}
\usepackage{amsfonts}
\usepackage{ifthen}
\usepackage{amsthm}
\usepackage{amsmath}
\usepackage{amssymb}
\usepackage{graphicx}
\usepackage{amscd,amssymb,amsthm}

\newcounter{minutes}
\divide\time by 60
\newcounter{hours}
\multiply\time by 60 \addtocounter{minutes}{-\time}

\title{Bisection of Geodesic Segments in Hyperbolic Geometry}

\author{Matti Vuorinen}
\author{Gendi Wang}

\address{Department of Mathematics and Statistics, University of Turku, Turku 20014,
Finland} \email{vuorinen@utu.fi}
\address{Department of Mathematics and Statistics, University of Turku, Turku 20014,
Finland}\email{genwan@utu.fi}


\newtheorem{theorem}[equation]{Theorem}

\newtheorem{lemma}[equation]{Lemma}
\newtheorem{proposition}[equation]{Proposition}

\newtheorem{corollary}[equation]{Corollary}
\newtheorem{remark}[equation]{Remark}

\newtheorem{propertyp}[equation]{Projection property}

\newcommand{\beq}{\begin{equation}}
\newcommand{\eeq}{\end{equation}}
\newcommand{\BB}{\mathbb{B}^2}
\newcommand{\RR}{\mathbb{R}^2}
\newcommand{\UH}{\mathbb{H}^2}
\newcommand{\Rn}{\mathbb{R}^n}
\newcommand{\R}{\mathbb{R}}

\numberwithin{equation}{section}

\begin{document}








\def\thefootnote{}
\footnotetext{ \texttt{\tiny File:~\jobname .tex,
          printed: \number\year-\number\month-\number\day,
          \thehours.\ifnum\theminutes<10{0}\fi\theminutes}
} \makeatletter\def\thefootnote{\@arabic\c@footnote}\makeatother

\maketitle

\begin{abstract}
Given a pair of points in the hyperbolic half plane or the unit disk,
we provide a simple construction of the midpoint of the hyperbolic geodesic
segment joining the points.
\end{abstract}

\section{Introduction}

Classical Euclidean geometry studies, in particular, configurations that
can be constructed in terms of compass and ruler. Thus, for instance,
given an angle we can bisect it and given a circle and a point outside
the corresponding disk, we can construct a tangent line to the circle
passing through the given point. For these facts and a general treatment
of the classical Euclidean geometry, see \cite{cg}.

In hyperbolic geometry, the parallel postulate no longer holds:
given a line and a point outside it, we can draw infinitely many lines
through the point, not intersecting the given line \cite{b}. On the other hand,
we can measure distances using the hyperbolic metric and consider
hyperbolic length minimizing curves, geodesic segments, between a
prescribed pair of points and form polygons whose sides consist of
geodesic segments, etc. In fact, many results of plane trigonometry
have their counterparts in this context in new vein, see \cite{b}. Recall that the hyperbolic area of a triangle with angles $\alpha, \beta,\gamma$ equals $\pi-(\alpha+ \beta+\gamma),$ so the sum of angles of
a triangle is always $< \pi\,$. For the history of hyperbolic geometry, we refer to \cite{d, gr,l, m}.

Our aim here is to study some basic constructions in hyperbolic geometry. As far as we know, in hyperbolic geometry constructions have been studied much less than in the Euclidean case and we have not been able to find our results in the literature. Our main results are the
following two theorems.

\begin{theorem}\label{my11}
Given a pair of points in the upper half plane $\UH\,,$ one can bisect
the hyperbolic segment joining the points by a geometric construction.
\end{theorem}

\begin{theorem}\label{my12}
Theorem \ref{my11} is also valid for the case of the unit disk $\BB\,.$
\end{theorem}


Previously, construction problems in hyperbolic geometry have been studied in \cite{go}. The bisection problem has been studied in \cite[Construction 3.1]{go} and in \cite[2.9]{kv} by use of methods different from ours. For an interesting survey of hyperbolic type
geometries, the reader is referred to \cite{himps, t}.

It is a basic fact that the hyperbolic geometries of the half plane and of the unit disk are isometrically equivalent via M\"obius transformations. Thus it is natural to expect that a construction in
one of these cases leads to a construction in the other case. However, our methods of construction, based on Euclidean compass and ruler, are
not M\"obius invariant. Because of this reason, we must treat these two cases separately.

The hyperbolic geodesic segment joining $x,y\in \BB$ is a subarc $J[x,y]$ of a circle perpendicular to $\partial{\BB}$. Therefore, the hyperbolic midpoint $z$ is determined as the set $\{z\}=J[x,y]\cap[0,e]$ for some point $e\in\partial{\BB}$, where $[0,e]=\{te: 0\leq t\leq 1\}$. We will find the vector $e$ by five different methods that produce five different points on the line determined by the radius $[0,e]$. Therefore, a byproduct of our five methods to prove Theorem 1.2 is the conclusion that these five points are on the same line, a fact, that may be of independent interest. Some of these points are shown in Figure \ref{cor}.

\section{Hyperbolic geometry}

The group of M\"obius transformations in $\overline{\Rn}$ is generated by transformations of two types:

(1) reflections in the hyperplane $P(a,t)=\{x\in \Rn: x\cdot a=t\}\cup\{\infty\}$

$$f_1(x)=x-2(x\cdot a-t)\frac{a}{|a|^2},\,\, f_1(\infty)=\infty,$$
where $a\in\Rn\setminus\{0\}$ and $t\in\R$;

(2) inversions (reflections) in the sphere $S^{n-1}(a,r)=\{x\in \Rn: |x-a|=r\}$

$$f_2(x)=a+\frac{r^2(x-a)}{|x-a|^2},\,\,f_2(a)=\infty, f_2(\infty)=a,$$
where $a\in\Rn$ and $r>0$.
If $G\subset\overline{\Rn}$ we denote by $\mathcal{GM}(G)$ the group of all M\"obius transformations which map $G$ onto itself.
For $x\in \Rn\setminus\{0\}$, we denote by $x^*=x/{|x|^2}$ the inversion in the sphere $S^{n-1}(0,1)=S^{n-1}$.

By symmetry, it is clear that the circle $S^1(a,r_a)$ containing $x, y, x^*, y^*$ is orthogonal to $S^1$. By
\cite[Exercise 1.1.27]{k} and \cite[Lemma 2.2]{kv}, we have
\beq\label{orar}
a =i\frac{y(1+|x|^2)-x(1+|y|^2)}{2(x_2y_1-x_1y_2)}\,\,\,\, {\rm and}\, \,\,\,r_a=\frac {|x-y|\big|x|y|^2-y\big|}{2|y||x_1y_2-x_2y_1|},
\eeq
for  $x=(x_1,x_2), y=(y_1,y_2)\in \mathbb{R}^2\setminus\{0\}$ such that $0,x,y$ are noncollinear.

The chordal distance is defined by
\beq\label{q}
\left\{\begin{array}{ll}
q(x,y)=\frac{|x-y|}{\sqrt{1+|x|^2}\sqrt{1+|y|^2}},&\,\,\, x\,,y\neq\infty,\\
q(x,\infty)=\frac{1}{\sqrt{1+|x|^2}},&\,\,\, x\neq\infty,
\end{array}\right.
\eeq
for $x\,,y\in\overline{\Rn}$.

For an ordered quadruple $a,b,c,d$ of distinct points in $\overline{\Rn}$ we define the absolute ratio by
$$|a,b,c,d|=\frac{q(a,c)q(b,d)}{q(a,b)q(c,d)}.$$
It follows from (\ref{q}) that for distinct points $a,b,c,d\in \Rn$

\beq\label{crossratio}
|a,b,c,d|=\frac{|a-c||b-d|}{|a-b||c-d|}.
\eeq
The most important property of the absolute ratio is M\"obius invariance, see \cite[Theorem 3.2.7]{b}, i.e., if $f$ is a M\"obius transformation, then
$$|f(a),f(b),f(c),f(d)|=|a,b,c,d|,$$
for all distinct $a,b,c,d\in \overline{\Rn}$.

Let $G\subsetneq\Rn$($n\geq 2$) be a domain and $w: G\rightarrow(0,\infty)$ be a continuous function.  We define the weighted length of a rectifiable curve $\gamma\subset G$ by
$$\ell_w(\gamma)=\int_{\gamma}w(z)|dz|$$
and the weighted distance by
$$d_w(x,y)=\inf_{\gamma}\ell_w(\gamma),$$
where the infimum is taken over all rectifiable curves in $G$ joining $x$ and $y$. It is easy to see that $d_w$ defines a metric on $G$ and $(G,d_w)$ is a metric space. We say that a curve $\gamma: [0,1]\rightarrow G$ is a geodesic joining $\gamma(0)$ and $\gamma(1)$ if for all $t\in (0,1)$, we have
$$d_w(\gamma(0),\gamma(1))=d_w(\gamma(0),\gamma(t))+d_w(\gamma(t),\gamma(1)).$$
The hyperbolic distance in $\UH$ is defined by the weight function $w_{\UH}(x)=1/{x_2}$ and in $\BB$ by the weight function $w_{\BB}(x)=2/{(1-|x|^2)}$. By \cite[p.35]{b} we have
\beq\label{cosh}
\cosh\rho_{\UH}(x,y)=1+\frac{|x-y|^2}{2x_2y_2}
\eeq
for all $x,y\in \UH$, and by \cite[p.40]{b} we have
\beq\label{sinh}
\sinh\frac{\rho_{\BB}(x,y)}{2}=\frac{|x-y|}{\sqrt{1-|x|^2}\sqrt{1-|y|^2}}
\eeq
for all $x,y\in \BB$.

In order to write formulas (\ref{cosh}) and (\ref{sinh}) in another form, let $\in G\in\{\UH,\BB\}$, $x,y\in{G}$. Let $L$ be an arc of a circle perpendicular to $\partial G$ with $x,y\in L$ and let $\{x_*,y_*\}=L\cap\partial G$, the points being labelled so that $x_*, x, y, y_*$ occur in this order on $L$. Then by \cite[(7.26)]{b}
\beq\label{rho}
\rho_G(x,y)=\sup\{\log|a,x,y,b|:a,b\in\partial G\}=\log|x_*,x,y,y_*|.
\eeq
We will omit the subscript $G$ if it is clear from the context.
Hyperbolic distance is invariant under M\"obius transformations of $G$ onto $G'$ for $G\,,G'\in\{\BB\,,\UH\}$.

\begin{figure}[h]
\begin{minipage}[t]{0.45\linewidth}
\centering
\includegraphics[width=8cm]{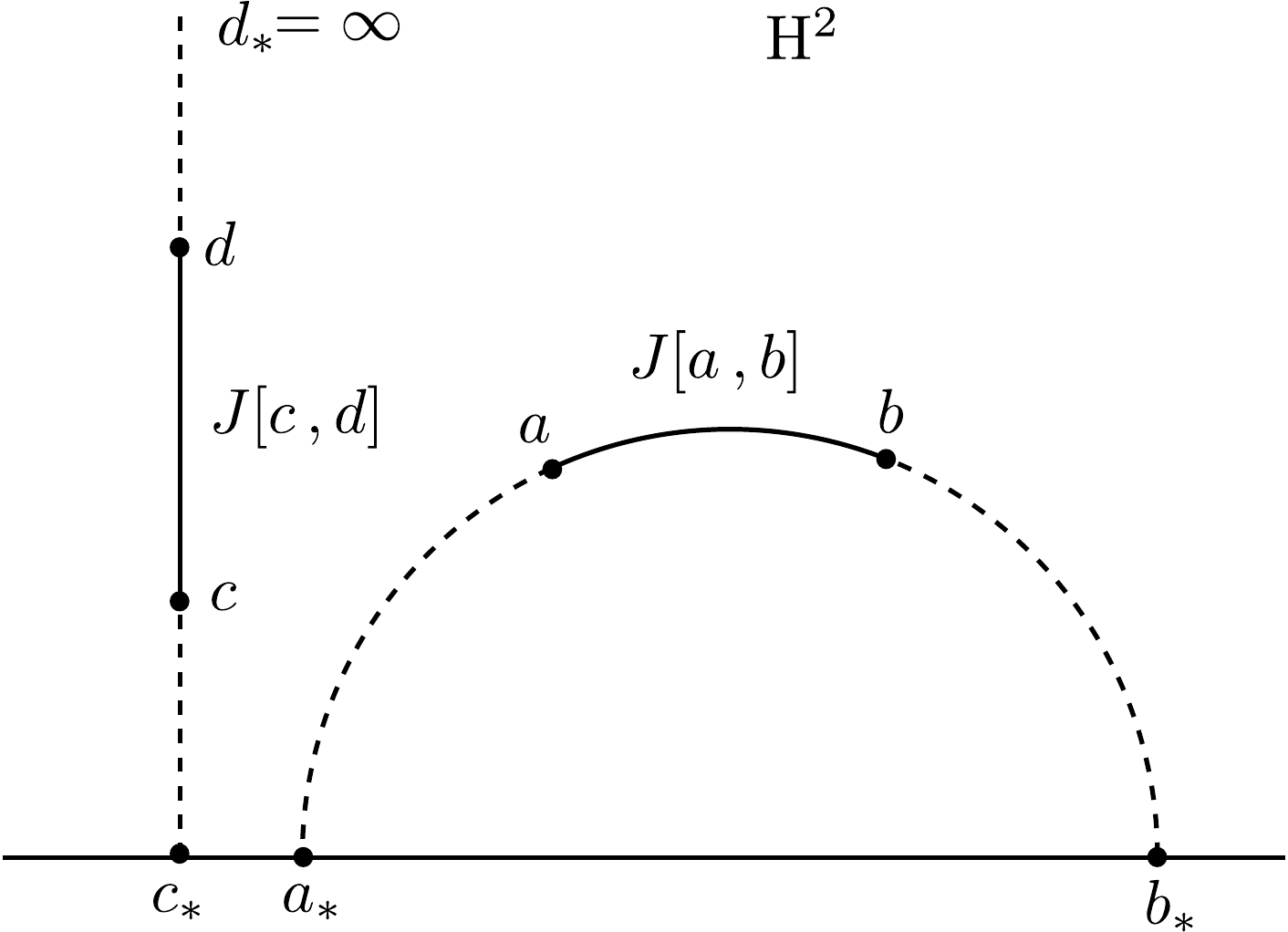}
\caption{\label{h2}}
\end{minipage}
\hfill
\hspace{1cm}
\begin{minipage}[t]{0.45\linewidth}
\centering
\includegraphics[width=7cm]{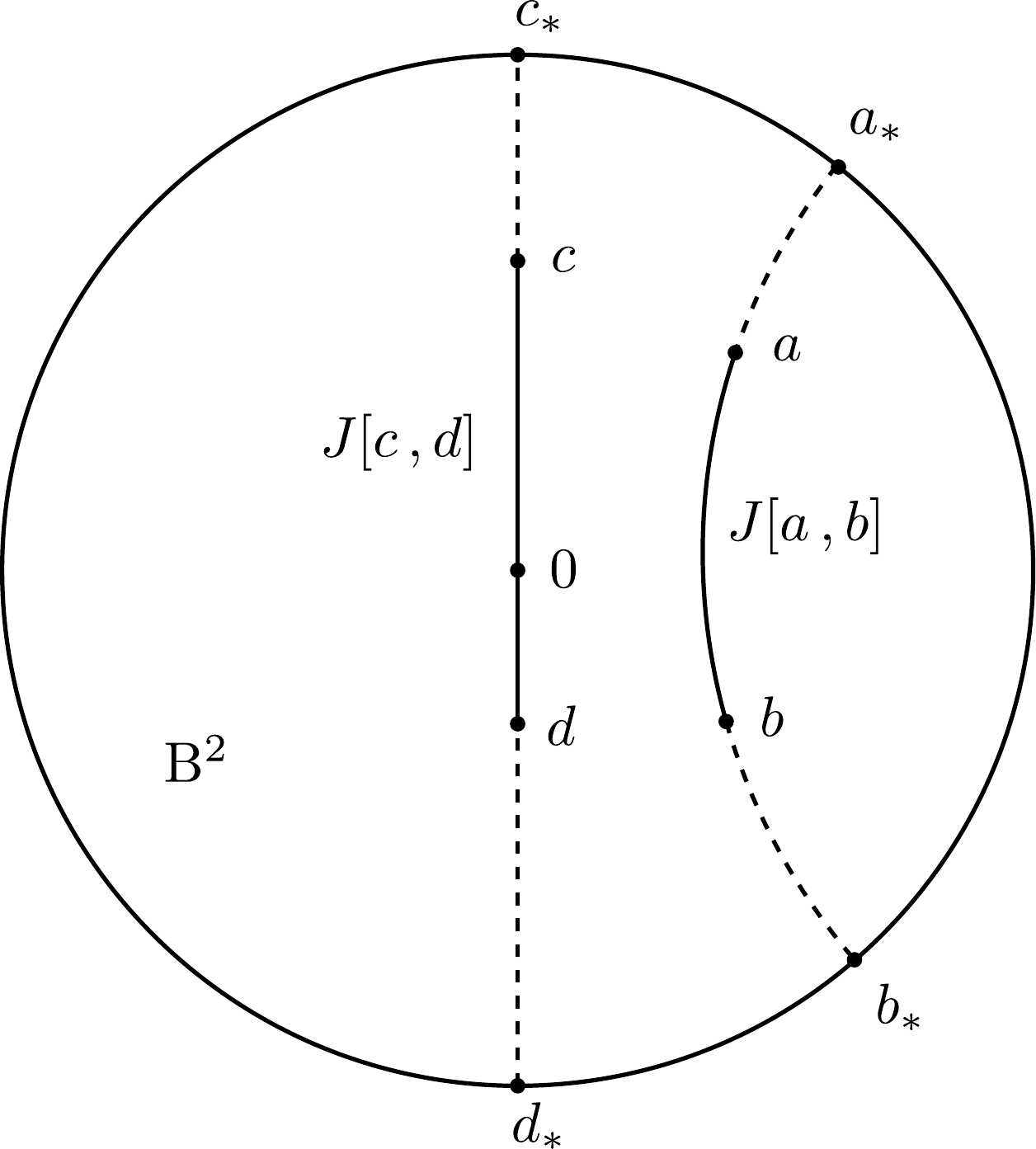}
\caption{ \label{b2}}
\end{minipage}
\end{figure}

Hyperbolic geodesics are  arcs of circles, which are orthogonal to the boundary of the domain. More precisely, for $a,b\in \BB$ (or $\UH)$, the hyperbolic geodesic segment joining $a$ and $b$ is an arc of a circle orthogonal to $S^1$ (or $\partial \UH)$. In a limiting case, the points $a$ and $b$ are located on a Euclidean line through $0$ (or located on a normal of $\partial \UH$), see \cite{b}.  Therefore, the points $x_*$ and $y_*$ are the end points of the hyperbolic geodesic. We denote by $J[a,b]$ the hyperbolic geodesic segment or shortly hyperbolic segment joining $a$ and $b$. For any two distinct points the hyperbolic geodesic segment is unique(see Figure \ref{h2} and \ref{b2}).

Knowing the geodesics, by \cite[p.21]{v} we have the hyperbolic distance  in two special cases in $\UH$.
First, for $r, s>0$, we have
\beq\label{logrs}
\rho(re_2, se_2)=\left|\int_s^r\frac{dt}{t}\right|=\left|\log \frac rs\right|.
\eeq
Second, if $\varphi\in(0,\pi/2)$ we denote $e_\varphi=(\cos \varphi) e_1+(\sin\varphi) e_2$ and have
\beq\label{loge}
\rho(e_2,e_\varphi)=\int_{J[e_\varphi,e_2]}\frac {d\alpha}{\sin \alpha}=\int_{\varphi}^\frac{\pi}{2}\frac {d\alpha}{\sin \alpha}=\log \mbox{cot}\frac12\varphi,
\eeq
where $e_1\,,e_2$ are the standard unit vectors in $\RR$. By (\ref{loge}), we calculate the midpoint $z$ of the hyperbolic segment $J[x,y]$ in $\UH$
\beq\label{midz}
z=e^{i\delta},\,\,\, \delta=\mbox{arc} \cos \left( \frac{\cos\frac{\beta+\alpha}{2}}{\cos\frac{\beta-\alpha}{2}}\right),
\eeq
where $x=e^{i\alpha}, y=e^{i\beta}$ and $0<\alpha<\beta<\pi$.

Next we give the counterparts of (\ref{logrs}) and (\ref{loge}) in $\BB$. By \cite[(2.17)]{v}, for $s\in(-t,t)$, we have
$$\rho(se_1,te_1)=\log \left(\frac{1+t}{1-t}\cdot\frac{1-s}{1+s}\right).$$

\begin{figure}[h]
\begin{minipage}[t]{0.45\linewidth}
\centering
\includegraphics[width=7.5cm]{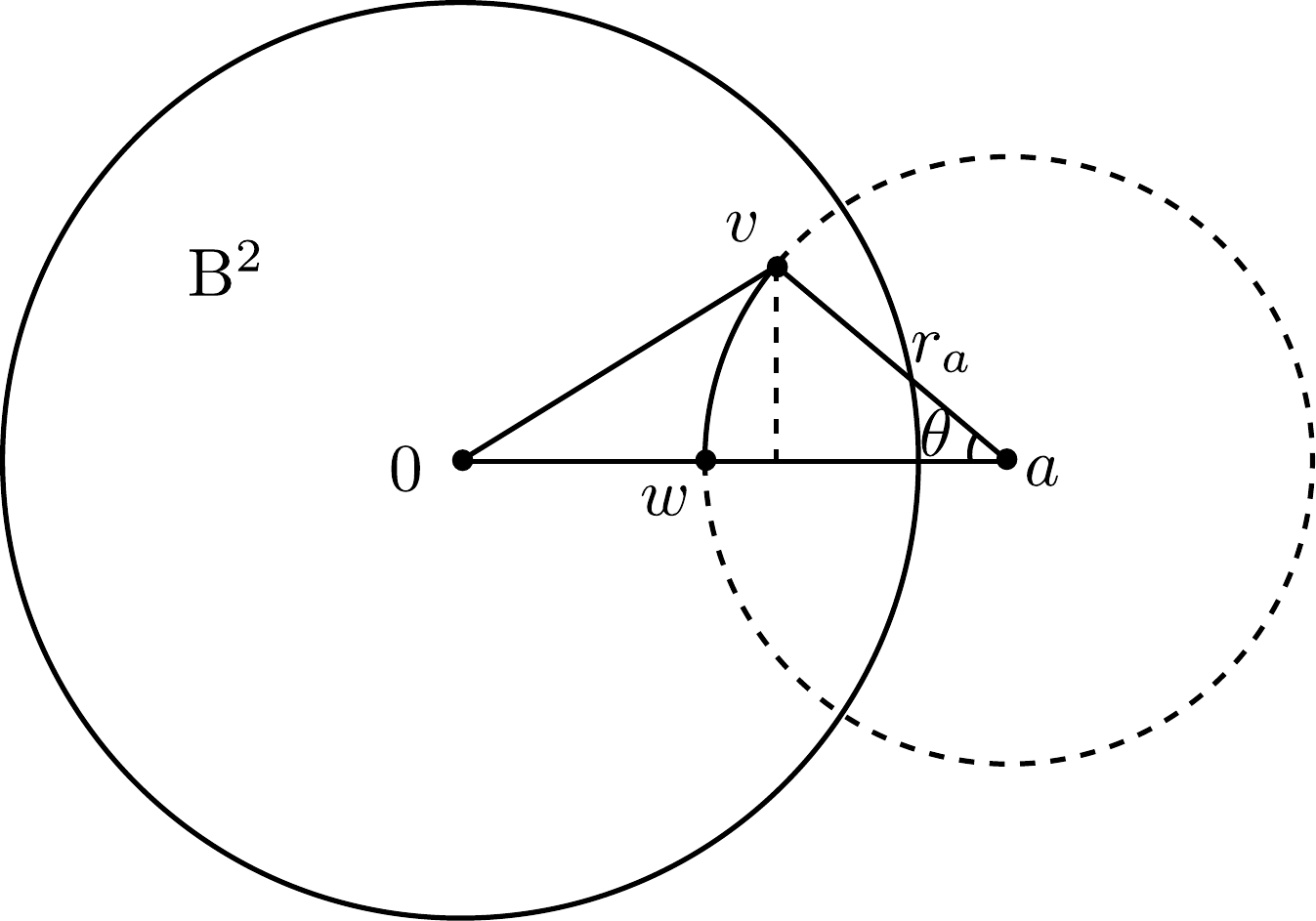}
\caption{ \label{vw}}
\end{minipage}
\end{figure}

For $v\in \BB\setminus\{0\}$ satisfying $0<{\rm arg} v<\pi/2$ and $0<w<1$, we have
\begin{eqnarray*}
|v|^2&=&(\sqrt{1+r^2_a}-r_a\cos\theta)^2+(r_a\sin\theta)^2\\
&=& 1+2r_a(r_a-\sqrt{1+r^2_a}\cos\theta),
\end{eqnarray*}
where $a,\,r_a$ are as in (\ref{orar}) by taking $x=v,\,y=w$ and $0<\theta=\angle 0av<\frac{\pi}{2}$, see Figure \ref{vw}.
Therefore, by the definition of the hyperbolic distance and \cite[4.3.133]{as}, we have
\begin{eqnarray}\label{rhowv}
\rho_{\BB}(w,v)&=&\int_0^{\theta}\frac{2r_adt}{2r_a(\sqrt{1+r^2_a}\cos t-r_a)}\nonumber\\
&=&\log\frac{A\tan\frac{\theta}{2}+1}{A\tan\frac{\theta}{2}-1},\,\,\,\,A=\sqrt{1+r^2_a}+r_a.
\end{eqnarray}
Let $x\,,y\in\BB\setminus\{0\}$ such that $-\pi/2<{\rm arg}x<{\rm arg} y<\pi/2$ and $0,x,y$ are noncollinear. Let $z$ be the midpoint of the hyperbolic segment $J[x,y]$. By (\ref{rhowv}), we have
$$\tan\frac{\delta}{2}=\frac{B+1}{A(B-1)}, \,\,\,\, B=\sqrt{\frac{(A\tan\frac{\beta}{2}+1)(A\tan\frac{\alpha}{2}-1)}{(A\tan\frac{\beta}{2}-1)(A\tan\frac{\alpha}{2}+1)}},$$
where $\alpha={\rm sgn}({\rm Im}\,x)\angle 0ax ,\beta={\rm sgn}({\rm Im}\,y)\angle 0ay,\delta={\rm sgn}({\rm Im}\,z)\angle 0az$, $\rm sgn$ is the signum function. It is easy to see that $-\frac{\pi}{2}<\alpha<\delta<\beta<\frac{\pi}{2}$.

\section{The half plane}
First we give some notation.

Let $L(x,y)$ be the line through the points $x$ and $y$, $L_{ab}(v)$ be the line through the point $v$ and orthogonal to the line $L(a,b)$. For simplicity of notation, we let $L(v)$ stand for the line through the point $v$ and orthogonal to $\partial \UH$. Let $\partial \UH$ be the real axis. For the convenience of proof, we exchange the complex number $x$ in two terms $x=x_1+ix_2$ and $x=r\,e^{i\,\alpha}$.

\begin{lemma}\label{le31}
Let $x,y\in S^1\cap \UH$ , $\{x_*,y_*\}=S^1\cap\partial \UH$, let $x_*,x,y,y_*$ occur in this order on $S^1$, and let $z$ be the midpoint of the hyperbolic segment $J[x,y]$, $w=L(x,y)\cap \partial \UH$, $v=L(x,x_*)\cap L(y,y_*)$. Then

(1) the line $L(w,z)$ is tangent to the circle $S^1$;

(2) the line $L(v,z)$ is orthogonal to $\partial \UH$;

(3) the line $L(a,z)$ is orthogonal to $\partial \UH$, where $a$ is as in (\ref{orar});

(4) $\angle y_1z_1y=\angle x_1z_1x$.
\end{lemma}

\begin{proof} Let $x=\cos \alpha+i\sin \alpha$, $y=\cos \beta+i\sin \beta$.

(1) By computation, we have
\beq\label{Lxy}
w_2-\sin \alpha=\frac{\sin\beta-\sin\alpha}{\cos\beta-\cos\alpha}(w_1-\cos \alpha).
\eeq
Putting $w_2=0$, we get
\beq\label{w1}
w_1=\frac{\sin(\beta-\alpha)}{\sin\beta-\sin\alpha}=\frac{\cos\frac{\beta-\alpha}{2}}{\cos\frac{\beta+\alpha}{2}}.
\eeq

Since
$$|z|^2+|w-z|^2=|w|^2\,\,\, \Leftrightarrow\,\,\, w\cdot z=1,$$
by (\ref{midz}) and (\ref{w1}), we get $w\cdot z=w_1z_1=1$. Therefore $L(w,z)$ is tangent to $S^1$.

(2) By symmetry we may assume that $\mbox{Re} x>\mbox{ Re} y$. Then $x_*=1$ and $y_*=-1$. It is easy to see that $v=v_1+iv_2$ satisfies the following equations:
\beq\label{V11V2}
\left\{\begin{array}{ll}
v_2=\frac{\sin\alpha}{\cos\alpha-1}(v_1-1)\\
v_2=\frac{\sin\beta}{\cos\beta+1}(v_1+1)
\end{array}\right..
\eeq
By (\ref{V11V2}), we have
\beq\label{v11}
v_1=\frac{\sin(\alpha+\beta)
+\sin \alpha-\sin\beta}{\sin(\alpha-\beta)
+\sin \alpha+\sin\beta}=\frac{\cos\frac{\beta+\alpha}{2}}{\cos\frac{\beta-\alpha}{2}}.
\eeq
Because $\mbox{Re} v=\mbox{Re} z$ by (\ref{midz}) and (\ref{v11}), we see that  $L(v,z)$ is orthogonal to $\partial \UH$.

(3) By (\ref{orar}), we have
\beq\label{a}
a=\frac{1}{\cos\frac{\beta-\alpha}{2}}e^{i\frac{\beta+\alpha}{2}}.
\eeq
Because $\mbox{Re} a=\mbox{Re} z$ by (\ref{midz}) and (\ref{a}), it is easy to see that $L(a,z)$ is orthogonal to $\partial \UH$.

(4) By (\ref{midz}), we have
\begin{eqnarray*}
& &\frac{y_2}{z_1-y_1}= \frac{x_2}{x_1-z_1}\\
&\Leftrightarrow& (x_2+y_2)z_1=x_1y_2+x_2y_1\\
&\Leftrightarrow&(\sin\alpha+\sin\beta)\frac{\cos\frac{\beta+\alpha}{2}}{\cos\frac{\beta-\alpha}{2}}=\sin(\alpha+\beta)\\
&\Leftrightarrow& 2\sin\frac{\beta+\alpha}{2}\cos\frac{\beta+\alpha}{2}=\sin(\alpha+\beta).
\end{eqnarray*}
Thus we conclude that $\angle y_1z_1y=\angle x_1z_1x$. This completes the proof.
\end{proof}

\begin{proposition}
Let $x,y,w,v$ be as in Lemma 3.1, $z=S^1\cap S^1(\frac{w}{2},\frac{|w|}{2})\cap \UH$ and $n=L(x,y)\cap S^1(\frac w2,\frac{|w|}{2})\cap \UH$. Let $\{s,t\}=S^1(a,r_a)\cap S^1(\frac{w}{2},\frac{|w|}{2})$, $u=L(a)\cap L(x,y)$, where $a$ and $r_a$ are as in (\ref{orar}). Then

(1) the point $z$ is the midpoint of the hyperbolic segment $J[x,y]$;

(2) the point $n$ is the midpoint of the Euclidean segment $[x,y]$;

(3) the point $v$ is on the circle $S^1(a,r_a)$;

(4) the circle $S^1(a,r_a)$ is orthogonal to the circle $S^1(\frac{w}{2},\frac{|w|}{2})$;

(5) the point $u$ is on the line $L(s,t)$.
\end{proposition}

\begin{proof} Let $x=\cos \alpha+i\sin \alpha$, $y=\cos \beta+i\sin \beta$. By symmetry we may assume that $0<\alpha<\beta<\pi$.

(1)It is easy to see that $z$ is the midpoint of $J[x,y]$ by Lemma 3.1(1).

(2)By similar triangles, we have
$$\frac{n_2}{n_1}=\frac{w_1-y_1}{y_2} \,\,\,\,,\,\,\,\,\frac{n_2}{y_2}=\frac{w_1-n_1}{w_1-y_1},$$
and hence
\beq\label{n1}
n_1=\frac{w_1y_2^2}{1+w_1^2-2w_1y_1}.
\eeq
By (\ref{w1}) and (\ref{n1}), we have
\begin{eqnarray*}
& & n_1=\frac{x_1+y_1}{2}\\
&\Leftrightarrow& (x_1+y_1)(w^2_1+1)-2(x_1y_1+1)w_1=0\\
&\Leftrightarrow& \cos^2\frac{\beta-\alpha}{2}+\cos^2\frac{\beta+\alpha}{2}-(\cos\alpha\cos\beta+1)=0\\
&\Leftrightarrow&\frac12(\cos(\beta-\alpha)+\cos(\beta+\alpha)+2)-(\cos\alpha\cos\beta+1)=0.
\end{eqnarray*}
Thus we obtain that $n$ is the midpoint of $[x,y]$.

(3) By Lemma \ref{le31}(2)(3), (\ref{V11V2})--(\ref{a}) and the orthogonality of $S^1(a,r_a)$ and $S^1$, we have
\begin{eqnarray*}
& &|v-a|^2=r^2_a\\
&\Leftrightarrow& v_2-a_2=\tan\frac{\beta-\alpha}{2}\\
&\Leftrightarrow& \sin\frac{\alpha-\beta}{2}+\sin\frac{\alpha+\beta}{2}-\sin\alpha\cos\frac{\beta-\alpha}{2}=\sin\frac{\beta-\alpha}{2}(\cos\alpha-1)\\
&\Leftrightarrow& \cos\alpha\sin\frac{\beta-\alpha}{2}+\sin\alpha\cos\frac{\beta-\alpha}{2}=\sin\frac{\alpha+\beta}{2}.
\end{eqnarray*}
Thus we obtain that $v$ is on $S^1(a,r_a)$.

(4)Since $r^2_a+(\frac{|w|}{2})^2=|a-\frac{w}{2}|^2\Leftrightarrow a\cdot w=1$ and the latter is true
by (\ref{w1}) and (\ref{a}), we obtain that $S^1(a,r_a)$ is orthogonal to $S^1(\frac{w}{2},\frac{|w|}{2})$.

(5) By the orthogonality of $S^1(a,r_a)$ and $S^1(\frac{w}{2},\frac{|w|}{2})$, it is easy to see that $L(s,t)$ is orthogonal to $L(a,\frac{w}{2})$. First observe that $m\in\{s,t\}$ satisfies the following equations:
$$\left\{\begin{array}{ll}
|m-a|^2=|a|^2-1\\
|m-\frac{w}{2}|^2=(\frac{|w|}{2})^2
\end{array}\right.,$$
then we have
\beq\label{m}
m\cdot(2a-w)=1.
\eeq
By (\ref{m}), $L(u,m)$ is orthogonal to $L(a,\frac{w}{2})$
\begin{eqnarray*}
&\Leftrightarrow& (u-m)\cdot(2a-w)=0\\
&\Leftrightarrow& u\cdot(2a-w)=1.
\end{eqnarray*}
By the assumption of $u$, (\ref{Lxy}) and (\ref{a}), we have
\beq\label{u1u2pr}
u_1=a_1\,\,\,\, {\rm and}\,\,\,\,u_2=\frac{1-a^2_1}{a_2}.
\eeq
Therefore $u\cdot(2a-w)=1$ by (\ref{w1}) and (\ref{u1u2pr}).
This completes the proof.
\end{proof}

\begin{remark}
The orthocenter $p$ of the triangle $\Delta{vx_*y_*}$ is on the circle $S^1(a,r_a)$, where $a$ and $r_a$ are as in (\ref{orar}).
\end{remark}

\noindent{\bf 3.12 Bisection of geodesic segment in $\UH$}

We now provide some constructions of the midpoint of the hyperbolic geodesic segment in $\UH$. Without loss of generality, we only need to deal with two cases: the points $x, y$ are on a line, which is orthogonal to $\partial \UH$ and the points $x, y$ are on a circle, which is orthogonal to $\partial \UH$.

{\bf Case 1.} {\it The hyperbolic segment $J[x,y]$ is located on the line $L(x,y)$ that is orthogonal to $\partial \UH$}

\begin{figure}[h]
\begin{minipage}[t]{0.45\linewidth}
\centering
\includegraphics[width=8cm]{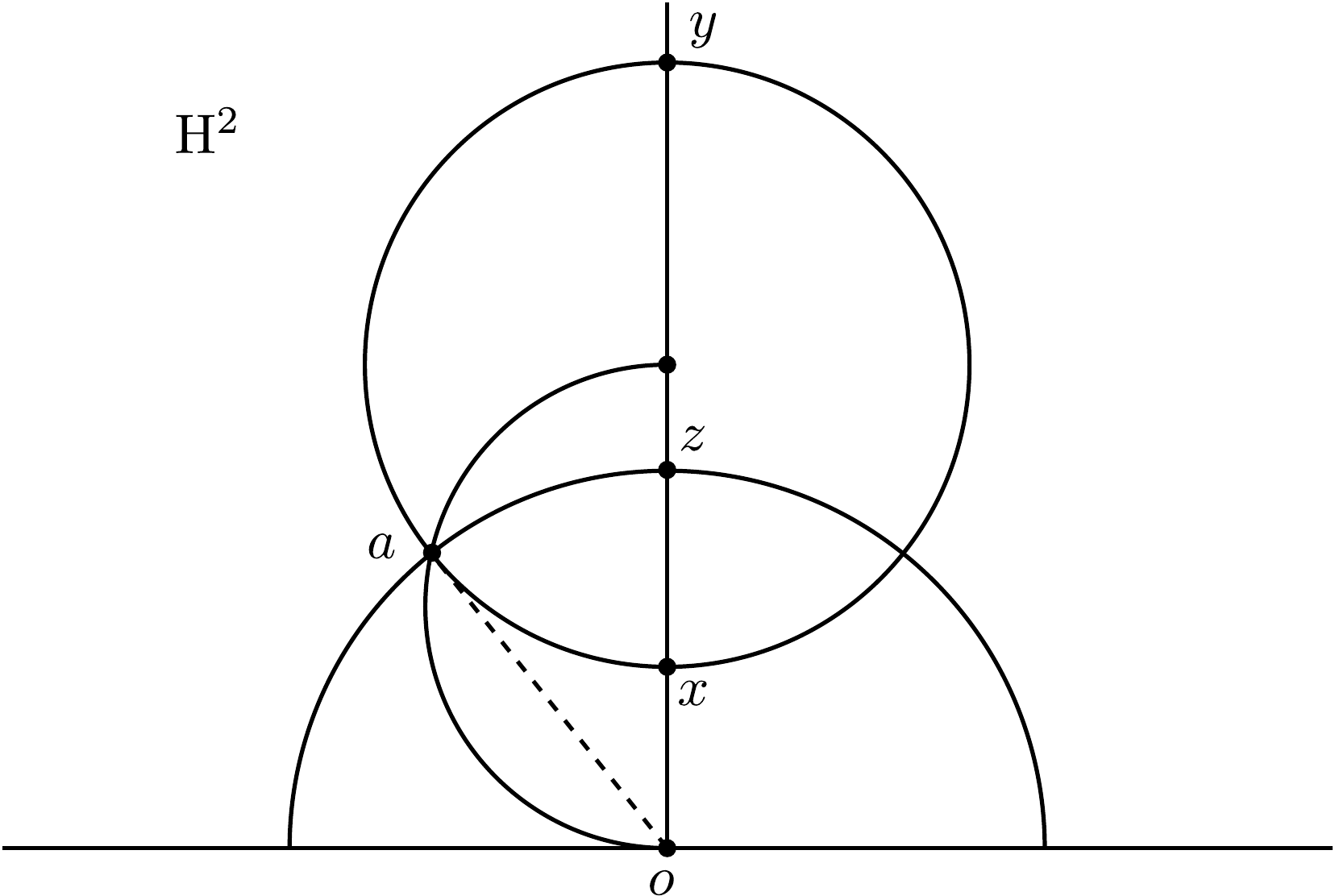}
\caption{ \label{31}}
\end{minipage}
\end{figure}

Without loss of generality, we may assume that ${\rm Im}\,x<{\rm Im}\,y$ and ${\rm Re}\,x={\rm Re}\,y=0$.

\vspace{-3.2mm}
\begin{table}[htbp]
\centering
\begin{tabularx}{\textwidth}{lX}
\emph{Step(1)} & Construct the line $L(x,y)$. Let $o=L(x,y)\cap \partial {\UH}$.\\
\emph{Step(2)} & Construct the circle $S^1(\frac{x+y}{2},\frac{|x-y|}{2})$.\\
\emph{Step(3)} & Construct the circle $S^1(\frac{x+y}{4},\frac{|x+y|}{4})$. Let $a=S^1(\frac{x+y}{4},\frac{|x+y|}{4})\cap S^1(\frac{x+y}{2},\frac{|x-y|}{2})$.\\
\emph{Step(4)} & Construct the circle $S^1(o,|a|)$.
\end{tabularx}
\end{table}
\vspace{-3.8mm}

\noindent By elementary geometry we see that $(\mbox{Im} z)^2=(\mbox{Im} x)(\mbox{Im} y)$ and therefore by (\ref{logrs}), we see that the midpoint $z$ of $J[x,y]$ is the intersection of $S^1(o,|a|)$ and $L(x,y)$ in $\UH$, see Figure \ref{31}.

{\bf Case 2.} {\it The hyperbolic segment $J[x,y]$ is located on the circle that is orthogonal to $\partial \UH$.}

{\noindent\bf Method I.}

\vspace{-3.2mm}
\begin{table}[htbp]
\centering
\begin{tabularx}{\textwidth}{lX}
\emph{Step (1)} & Construct the circle $S^1(o,r)$, which is orthogonal to $\partial \UH$ and contains the points $x,y$.\\
\emph{Step (2)} & Construct the line $L(x,y)$. Let $w=L(x,y)\cap \partial \UH$.\\
\emph{Step (3)} & Construct the circle $S^1(\frac{w+o}{2},\frac{|w-o|}{2})$.
\end{tabularx}
\end{table}
\vspace{-3.8mm}

\begin{figure}[h]
\begin{minipage}[t]{0.45\linewidth}
\centering
\includegraphics[width=8.2cm]{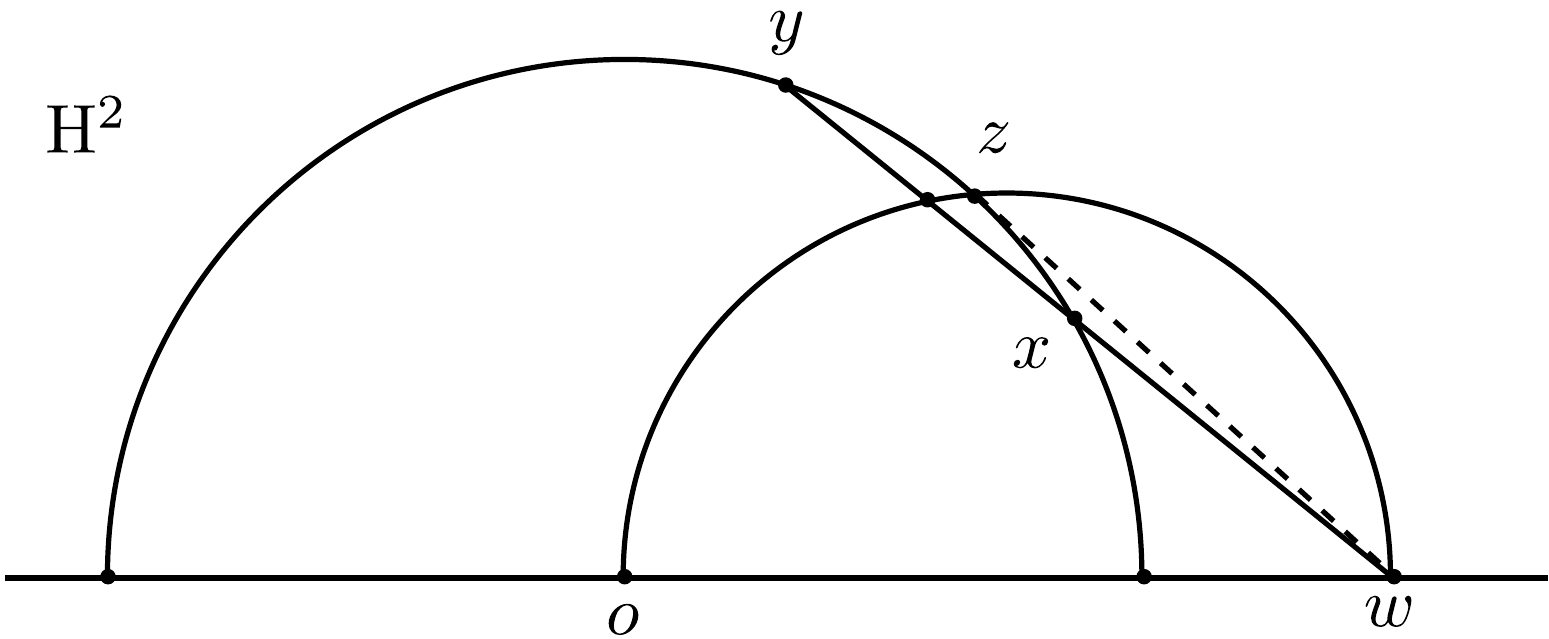}
\caption{\label{32}}
\end{minipage}
\hfill
\begin{minipage}[t]{0.45\linewidth}
\centering
\includegraphics[width=8cm]{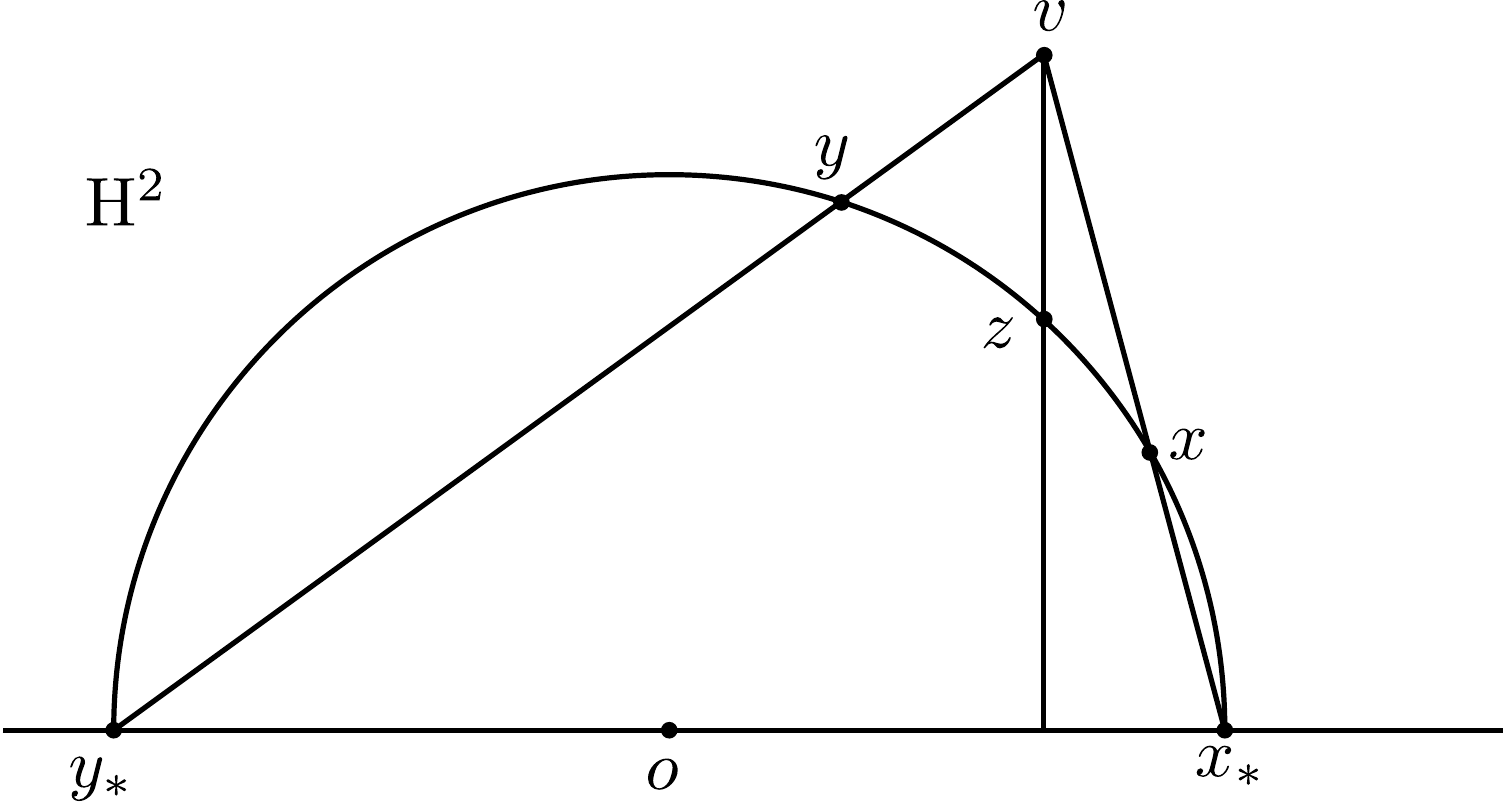}
\caption{\label{33}}
\end{minipage}
\end{figure}

\noindent Then the midpoint $z$ of $J[x,y]$ is the intersection of $S^1(\frac{w+o}{2},\frac{|w-o|}{2})$ and $S^1(o,r)$ in $\UH$ by Lemma 3.1(1), see Figure \ref{32}.

{\noindent\bf Method II.}

\vspace{-3.2mm}
\begin{table}[htbp]
\centering
\begin{tabularx}{\textwidth}{lX}
\emph{Step (1)} & Construct the circle $S^1(o,r)$, which is orthogonal to $\partial \UH$ and contains the points $x,y$.
Let $\{x_*,y_*\}=S^1(o,r)\cap \partial \UH$, $x_*,x,y,y_*$ occur in this order on $S^1(o,r)$.\\
\emph{Step (2)} & Construct the lines $L(x,x_*)$ and $L(y,y_*)$. Let $v=L(x,x_*)\cap L(y,y_*)$.\\
\emph{Step (3)} & Construct the line $L(v)$.
\end{tabularx}
\end{table}
\vspace{-3.8mm}

\noindent Then the midpoint $z$ of $J[x,y]$ is the intersection  of $L(v)$ and $S^1(o,r)$ in $\UH$ by Lemma 3.1(2), see Figure \ref{33}.

{\noindent\bf Method III.}

\vspace{-3.2mm}
\begin{longtable}{ll}
\centering
\emph{Step (1)} & Construct\, the\, circle\, $S^1(o,r)$,\, which\, is\, orthogonal\, to\, $\partial \UH$\, and\, contains\, the\\
                & points $x,y$.\\
\emph{Step (2)} & Construct the circle $S^1(a,r_a)$, which is orthogonal to $S^1(o,r)$ and contains the\\
                & points $x,y$.\\
\emph{Step (3)} & Construct the line $L(a)$.
\end{longtable}

\newpage

\begin{figure}[h]
\begin{minipage}[t]{0.45\linewidth}
\centering
\includegraphics[width=7.8cm]{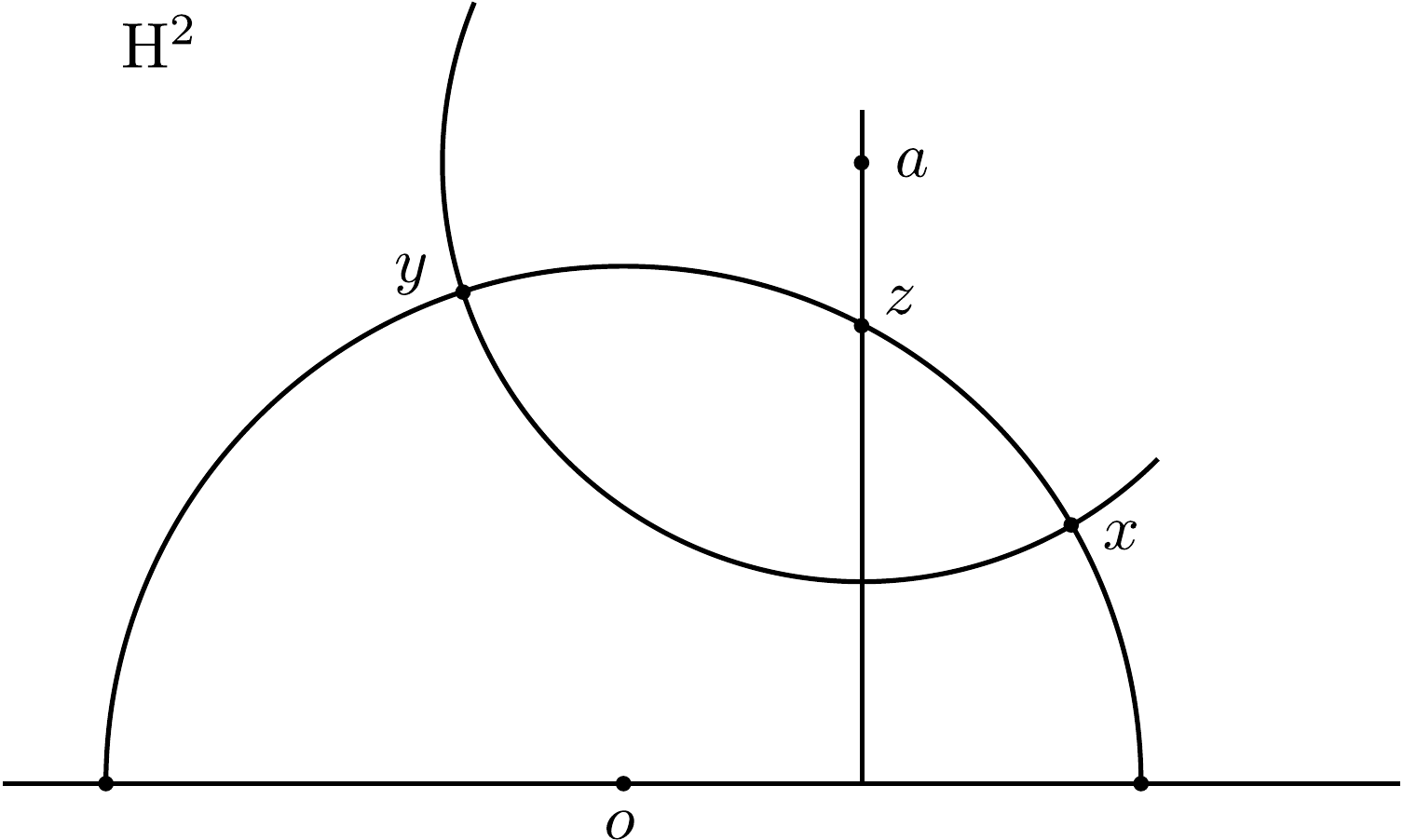}
\caption{\label{34}}
\end{minipage}
\hfill
\begin{minipage}[t]{0.45\linewidth}
\centering
\includegraphics[width=7cm]{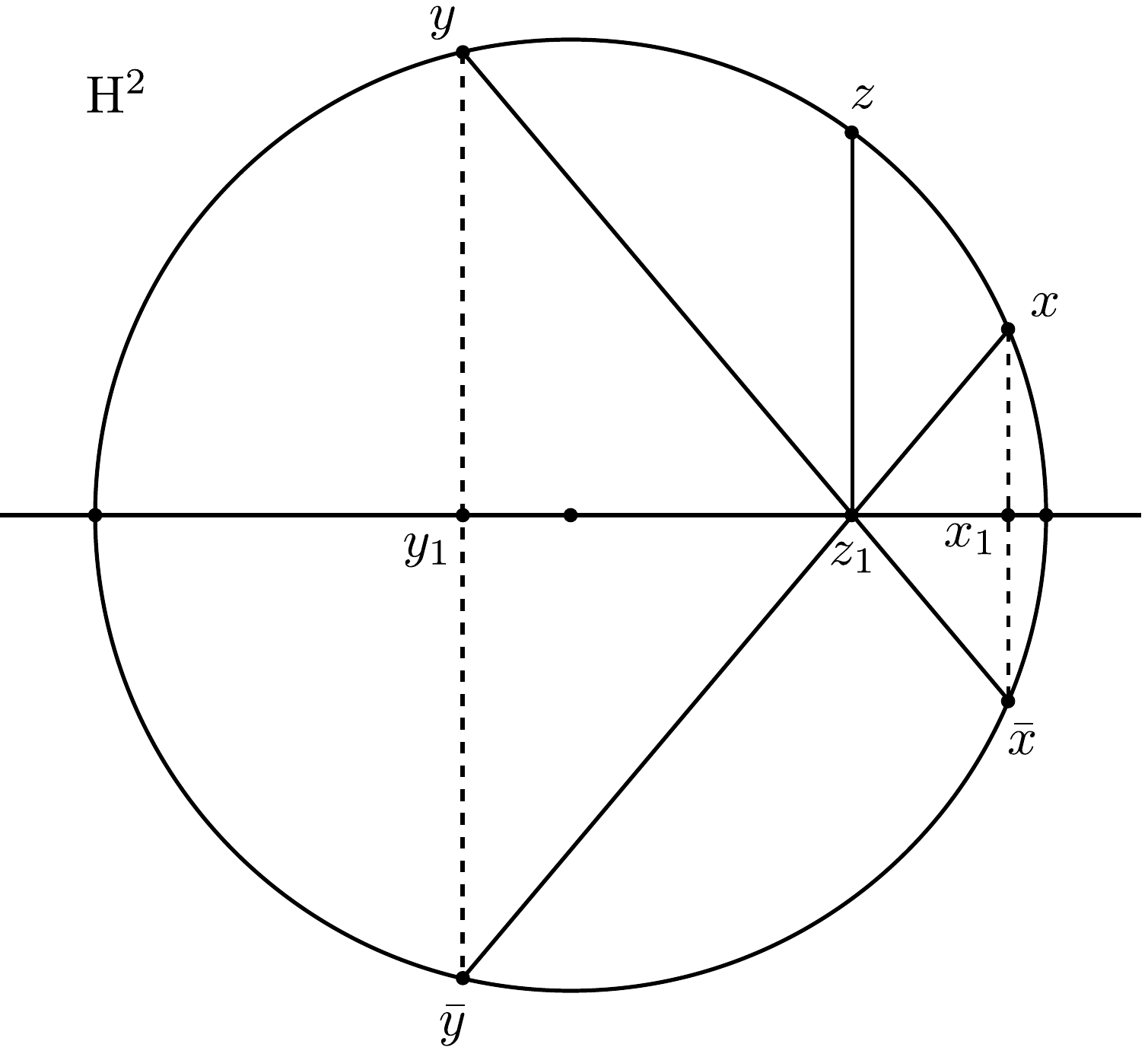}
\caption{\label{35}}
\end{minipage}
\end{figure}

\noindent Then the midpoint $z$ of $J[x,y]$ is the intersection  of $L(a)$ and $S^1(o,r)$ in $\UH$ by Lemma 3.1(3), see Figure \ref{34}.

{\noindent\bf Method IV.}

\vspace{-3.2mm}
\begin{table}[htbp]
\centering
\begin{tabularx}{\textwidth}{lX}
\emph{Step (1)} & Construct the circle $S^1(o,r)$, which is orthogonal to $\partial \UH$ and contains the points $x,y$.\\
\emph{Step (2}) & Construct the lines $L(x,\overline{y})$ and $L(\overline{x},y)$. Let $z_1=L(x,\overline{y})\cap L(\overline{x},y)\cap \partial \UH$.\\
\emph{Step (3)} & Construct the line $L(z_1)$.
\end{tabularx}
\end{table}
\vspace{-3.8mm}

\noindent Then the midpoint $z$ of $J[x,y]$ is the intersection  of $L(z_1)$ and $S^1(o,r)$ in $\UH$ by Lemma 3.1(4), see Figure \ref{35}.\\

\section{The unit disk}

We first show the relation between $\rho_{\UH}$ and $\rho_{\BB}$.

\begin{propertyp}
For $x,y\in S^1\cap \UH$,
\beq\label{rhoHB}
2\rho_{\UH}(x,y)=\rho_{\BB}(Pr(x),Pr(y)),
\eeq
where $Pr(x)$ and $Pr(y)$ are the points of the projection of  $x$ and $y$ on $\partial \UH$, respectively.
\end{propertyp}

\begin{proof} By symmetry we may assume that ${\rm Re}\,x<{\rm Re}\,y$. Then by similar triangles and (\ref{crossratio}), we have
$$|-1,x,y,1|^2=|-1,Pr(x),Pr(y),1|,$$
thus we
obtain (\ref{rhoHB}) by (\ref{rho}).
\end{proof}

\begin{remark} By the above property, we obtain that the projection of the midpoint $z$ of the hyperbolic segment $J[x,y]$ on the circle $S^1\cap\ \UH$ on $\partial \UH$ is just the midpoint $z_1$ of the hyperbolic segment $J[x_1,y_1]$ in $\BB$, which is also the projection of $J[x,y]$ on $\partial \UH$(cf. Figure \ref{35}).
\end{remark}

\begin{corollary} Given $X_1\in \BB\setminus\{0\}$, we can construct a sequence of points $\{X_k\}$ on the same radius with
\beq\label{kc}
\rho_{\BB}(0,X_k)=kc,
\eeq
where $c=\rho_{\BB}(0,\,X_1)$.
\end{corollary}

\begin{proof}
\emph{Step(1)} Construct the lines $L(0,X_1)$, $L_{0X_1}(0)$ and $L_{0X_1}(X_1)$.
Let $\{M_0,N_0\}=L_{0X_1}(0)\cap S^1$, $\{M_1,N_1\}=L_{0X_1}(X_1)\cap S^1$.

\emph{Step(2)} Construct the line $L(M_0,X_1)$.
Let $N_2=L(M_0,X_1)\cap S^1$.

\emph{Step(3)} Construct the line $L_{0X_1}(N_2)$. Let $M_2=L_{0X_1}(N_2)\cap S^1$, $X_2=L_{0X_1}(N_2)\cap L(0,X_1)$.

Repeat Steps (2) and (3), then we get a sequence of points $\{X_k\}$, which satisfy (\ref{kc}).
\end{proof}

\begin{lemma}
Let $x,y\in \BB\setminus\{0\}$ such that $0,x,y$ are noncollinear and $|x|\neq|y|$. Let $\{x_*,y_*\}=S^1\cap S^1(a,r_a)$, where $a,r_a$ are as in (2.1), then $x_*,x,y,y_*$ occur in this order on the orthogonal circle $S^1(a,r_a)$. Let $w=L(x,y)\cap L(x^*,y^*)$, $u=L(x,y^*)\cap L(y,x^*)$, $v=L(x,x_*)\cap L(y,y_*)$, $s=L(x,y_*)\cap L(y,x_*)$, $t=L(x_*,y^*)\cap L(y_*,x^*)$, $k=L(x_*,x^*)\cap L(y_*,y^*)$ . Construct the circle $S^1(w,r_w)$, which is orthogonal to the circle $S^1(a,r_a)$, let $z=S^1(w,r_w)\cap S^1(a,r_a)\cap \BB$. Then

(1) the circle $S^1(w,r_w)$ is orthogonal to the circle $S^1$ and
the point $z$ is the midpoint of the hyperbolic segment $J[x,y]$, where
\beq\label{w}
w=\frac{y(1-|x|^2)-x(1-|y|^2)}{|y|^2-|x|^2}\,\,\,\, and\,\,\,\, r_w=\frac{|x-y|\sqrt{(1-|x|^2)(1-|y|^2)}}{\big||y|^2-|x|^2\big|};
\eeq

(2) the points $v,s,t,k$ are on the same line $L(0,z)$ and $u$ is the point of the intersection of the lines $L(0,z)$ and $L(x_*,y_*)$, where
\beq\label{u}
u=\frac{y(1-|x|^2)+x(1-|y|^2)}{1-|x|^2|y|^2}.
\eeq
\end{lemma}

\begin{proof} (1) By computation, it is easy to see that $w=w_1+iw_2$ satisfies the following equations:
\beq\label{W1W2}
\left\{\begin{array}{ll}
w_2-y_2=\frac{y_2-x_2}{y_1-x_1}(w_1-y_1)\\
w_2-\frac{y_2}{|y|^2}=\frac{y_2|x|^2-x_2|y|^2}{y_1|x|^2-x_1|y|^2}(w_1-\frac{y_1}{|y|^2})
\end{array}\right..
\eeq
Solving (\ref{W1W2}), we have
\beq\label{w1w2}
w_1=\frac{y_1(1-|x|^2)-x_1(1-|y|^2)}{|y|^2-|x|^2}\,\,\, \,{\rm and} \,\,\,\,
w_2=\frac{y_2(1-|x|^2)-x_2(1-|y|^2)}{|y|^2-|x|^2},
\eeq
thus we obtain $w$ by (\ref{w1w2}).

Since
\begin{eqnarray*}
& &r^2_w+1=|w|^2\\
&\Leftrightarrow& |w-a|^2-(|a|^2-1)+1=|w|^2\\
&\Leftrightarrow& w\cdot a=1,
\end{eqnarray*}
by (\ref{orar}) and (\ref{w1w2}), we have $w\cdot a=1$. Thus we prove that $S^1(w,r_w)$ is orthogonal to $S^1$ and hence we get $r_w$. By \cite[Corollary 2.7]{kv} and (\ref{w1w2}), we see that $z$ is the midpoint of $J[x,y]$.

(2) Without loss of generality, we may assume that $a_2=0$. Then we have
\beq\label{|x||y|}
1-|x|^2=2-\frac{x_1}{y_1-x_1}(|y|^2-|x|^2)\,\,\,\, {\rm and} \,\,\,\, 1-|y|^2=2-\frac{y_1}{y_1-x_1}(|y|^2-|x|^2).
\eeq
By (\ref{w1w2}) and (\ref{|x||y|}), we also have
\beq\label{w1w2a2=0}
w_1=\frac{2(y_1-x_1)}{|y|^2-|x|^2}\,\,\,\, {\rm and} \,\,\,\,w_2=\frac{2(y_2-x_2)}{|y|^2-|x|^2}-\frac{x_1y_2-x_2y_1}{y_1-x_1}.
\eeq

By \cite[1.34]{v}, there exists a sense-preserving M\"obius transformation $T$ such that $T(\BB)=\BB$ and $T(z)=0$. By the orthogonality of $S^1(a,r_a)$ and $S^1$, we have $T(x)=-T(y)$ and consequently $T(S^1(a,r_a))=L(T(x),T(y))$. Since $S^1(w,r_w)$ is orthogonal to both $S^1(a,r_a)$ and $S^1$, we conclude that $T(S^1(w,r_w))=L_{T(x)T(y)}(T(z))$. Therefore, it is easily seen that $T(z')=\infty$, where $z'$ is another point of the intersection of $S^1(a,r_a)$ and $S^1(w,r_w)$. Then we have that the line containing $T(0),T(z),T(z')$ is orthogonal to $T(S^1)=S^1$. Because M\"obius transformations preserve circles and angles, we obtain that $0,z,z'$ are collinear and hence $L(0,z)$ is orthogonal to $L(w,a)$.
The same reason gives that $w,x_*,y_*$ are collinear, so we may assume that $x_*=w_1+im$ and $y_*=w_1-im$, $m>0$. We also have
\beq\label{ma1}
m^2=1-w^2_1 \,\,\,\,{\rm and}\,\,\,\, a_1=\frac{1}{w_1}.
\eeq

First, we prove that $v$ is on $L(0,z)$.
We observe that $v=v_1+iv_2$ satisfies the following equations:
\beq\label{V1V2}
\left\{\begin{array}{ll}
v_2-m=\frac{x_2-m}{x_1-w_1}(v_1-w_1)\\
v_2+m=\frac{y_2+m}{y_1-w_1}(v_1-w_1)
\end{array}\right..
\eeq
Solving (\ref{V1V2}), we have
\beq\label{v1}
v_1=w_1+2m\frac{(x_1-w_1)(y_1-w_1)}{(x_1-w_1)(y_2+m)-(x_2-m)(y_1-w_1)}
\eeq
and
\beq\label{v2}
v_2=m\frac{(x_1-w_1)(y_2+m)+(x_2-m)(y_1-w_1)}{(x_1-w_1)(y_2+m)-(x_2-m)(y_1-w_1)}.
\eeq
By (\ref{ma1}), we have $L(0,v)$ is orthogonal to $L(w,a)$
\begin{eqnarray*}
&\Leftrightarrow& (w-a)\cdot v=0\\
&\Leftrightarrow& m^2v_1-w_1w_2v_2=0\\
&\Leftrightarrow& 2m^2(x_1-w_1)(y_1-w_1)+w_1\{(x_1+y_1-2w_1)m^2+[(x_1-w_1)(y_2-w_2)\\
&   &-(x_2-w_2)(y_1-w_1)]m-w_2[y_2(x_1-w_1)+x_2(y_1-w_1)]\}=0.
\end{eqnarray*}
By similar triangles, we have $(x_1-w_1)(y_2-w_2)=(x_2-w_2)(y_1-w_1)$. The task is now to prove
\beq\label{uvst}
(1-w^2_1)[2x_1y_1-(x_1+y_1)w_1]+w_1w_2[(x_2+y_2)w_1-(x_1y_2+x_2y_1)]=0.
\eeq

By (\ref{w1w2a2=0}) and simplification, we get
\begin{eqnarray}\label{416}
& &{\rm Eq.} (\ref{uvst})\nonumber\\
&\Leftrightarrow&[(|y|^2-|x|^2)^2-4(y_1-x_1)^2][x_1y_1(|y|^2-|x|^2)-(y^2_1-x^2_1)]\\
& &+[(x_1y_2+x_2y_1)(|y|^2-|x|^2)-2(y_1-x_1)(y_2+x_2)]\nonumber\\
& &\times[(x_1y_2-x_2y_1)(|y|^2-|x|^2)-2(y_1-x_1)(y_2-x_2)]=0\nonumber\\
&\Leftrightarrow& (|y|^2-|x|^2)^2D+4(|y|^2-|x|^2)(y_1-x_1)E=0,\nonumber
\end{eqnarray}
where
$$D=x_1y_1(|y|^2-|x|^2)+x^2_1-y^2_1+x^2_1y^2_2-x^2_2y^2_1$$
and
$$E=(y_1-x_1)(1-x_1y_1)-(x_1y^2_2-x^2_2y_1).$$
Then $D=(x_1+y_1)[x_1(1+|y|^2)-y_1(1+|x|^2)]=0$ and $E=y_1(1+|x|^2)-x_1(1+|y|^2)=0$ by the assumption of $a_2=0$. Hence, we obtain that $v$ is on $L(0,z)$.

\medskip

The proof for the result involving $s$ is almost the same as that for $v$.

\medskip

Next, we prove that $t$ is on $L(0,z)$.
We observe that $t=t_1+it_2$ satisfies the following equations:
\beq\label{T1T2}
\left\{\begin{array}{ll}
t_2+m=\frac{x_2+m|x|^2}{x_1-w_1|x|^2}(t_1-w_1)\\
t_2-m=\frac{y_2-m|y|^2}{y_1-w_1|y|^2}(t_1-w_1)
\end{array}\right..
\eeq
Solving (\ref{T1T2}), we have
\beq\label{t1}
t_1=w_1-2m\frac{(x_1-w_1|x|^2)(y_1-w_1|y|^2)}{(x_1-w_1|x|^2)(y_2-m|y|^2)-(x_2+m|x|^2)(y_1-w_1|y|^2)}
\eeq
and
\beq\label{t2}
t_2=-m\frac{(x_1-w_1|x|^2)(y_2-m|y|^2)+(x_2+m|x|^2)(y_1-w_1|y|^2)}{(x_1-w_1|x|^2)(y_2-m|y|^2)-(x_2+m|x|^2)(y_1-w_1|y|^2)}.
\eeq
By (\ref{w1w2})--(\ref{w1w2a2=0}) and the assumption of $a_2=0$, we have
\beq\label{x1-w1|y|^2}
x_1-w_1|x|^2=\frac{(y_1-x_1)(1-|x|^2)}{|y|^2-|x|^2}=w_1-x_1
\eeq
and
\beq\label{y1-w1|y|^2}
y_1-w_1|y|^2=\frac{(y_1-x_1)(1-|y|^2)}{|y|^2-|x|^2}=w_1-y_1.
\eeq
By (\ref{w1w2}), we also have
\beq\label{x2-w2|x|^2}
x_2-w_2|x|^2=\frac{(x_2|y|^2-y_2|x|^2)(1-|x|^2)}{|y|^2-|x|^2}
\eeq
and
\beq\label{y2-w2|y|^2}
y_2-w_2|y|^2=\frac{(x_2|y|^2-y_2|x|^2)(1-|y|^2)}{|y|^2-|x|^2}.
\eeq

By (\ref{ma1}) and (\ref{x1-w1|y|^2})--(\ref{y2-w2|y|^2}) , we have that $L(0,t)$ is orthogonal to $L(w,a)$
\begin{eqnarray*}
&\Leftrightarrow& (w-a)\cdot t=0\\
&\Leftrightarrow& m^2t_1-w_1w_2t_2=0\\
&\Leftrightarrow& m^2[y_1(x_1-w_1)+x_1(y_1-w_1)]-mw_1[(x_1-w_1)(y_2-w_2|y|^2)-\\
& &(y_1-w_1)(x_2-w_2|x|^2)]-w_1w_2[(x_1-w_1)y_2+x_2(y_1-w_1)]=0\\
&\Leftrightarrow& (\ref{uvst}).
\end{eqnarray*}
Thus we obtain that $t$ is on $L(0,z)$.

\medskip

The proof for the result involving $k$ is almost the same as that for $t$.

\medskip

Finally, we prove the result for $u$.
We observe that $u=u_1+iu_2$ satisfies the following equations:
\beq\label{U1U2}
\left\{\begin{array}{ll}
u_2-x_2=\frac{y_2-x_2|y|^2}{y_1-x_1|y|^2}(u_1-x_1)\\
u_2-y_2=\frac{x_2-y_2|x|^2}{x_1-y_1|x|^2}(u_1-y_1)
\end{array}\right..
\eeq
Solving (\ref{U1U2}), we have
\beq\label{u1}
u_1=\frac{y_1(1-|x|^2)+x_1(1-|y|^2)}{1-|x|^2|y|^2}
\eeq
and
\beq\label{u2}
u_2=\frac{y_2(1-|x|^2)+x_2(1-|y|^2)}{1-|x|^2|y|^2}.
\eeq
Hence we get (\ref{u}) by (\ref{u1}) and (\ref{u2}).

By (\ref{w1w2a2=0}), (\ref{u1}) and (\ref{u2}), we have that $L(0,u)$ is orthogonal to $L(w,a)$
\begin{eqnarray*}
&\Leftrightarrow& (w-a)\cdot u=0\\
&\Leftrightarrow& m^2u_1-w_1w_2u_2=0\\
&\Leftrightarrow& (\ref{416}).
\end{eqnarray*}
By (\ref{w1w2a2=0}), (\ref{u1}) and the assumption of $a_2=0$, we have $u_1=w_1$. Thus $u$ is the intersection of $L(x_*,y_*)$ and $L(0,z)$.
This completes the proof.
\end{proof}

\begin{figure}[h]
\centering
\includegraphics[width=14cm]{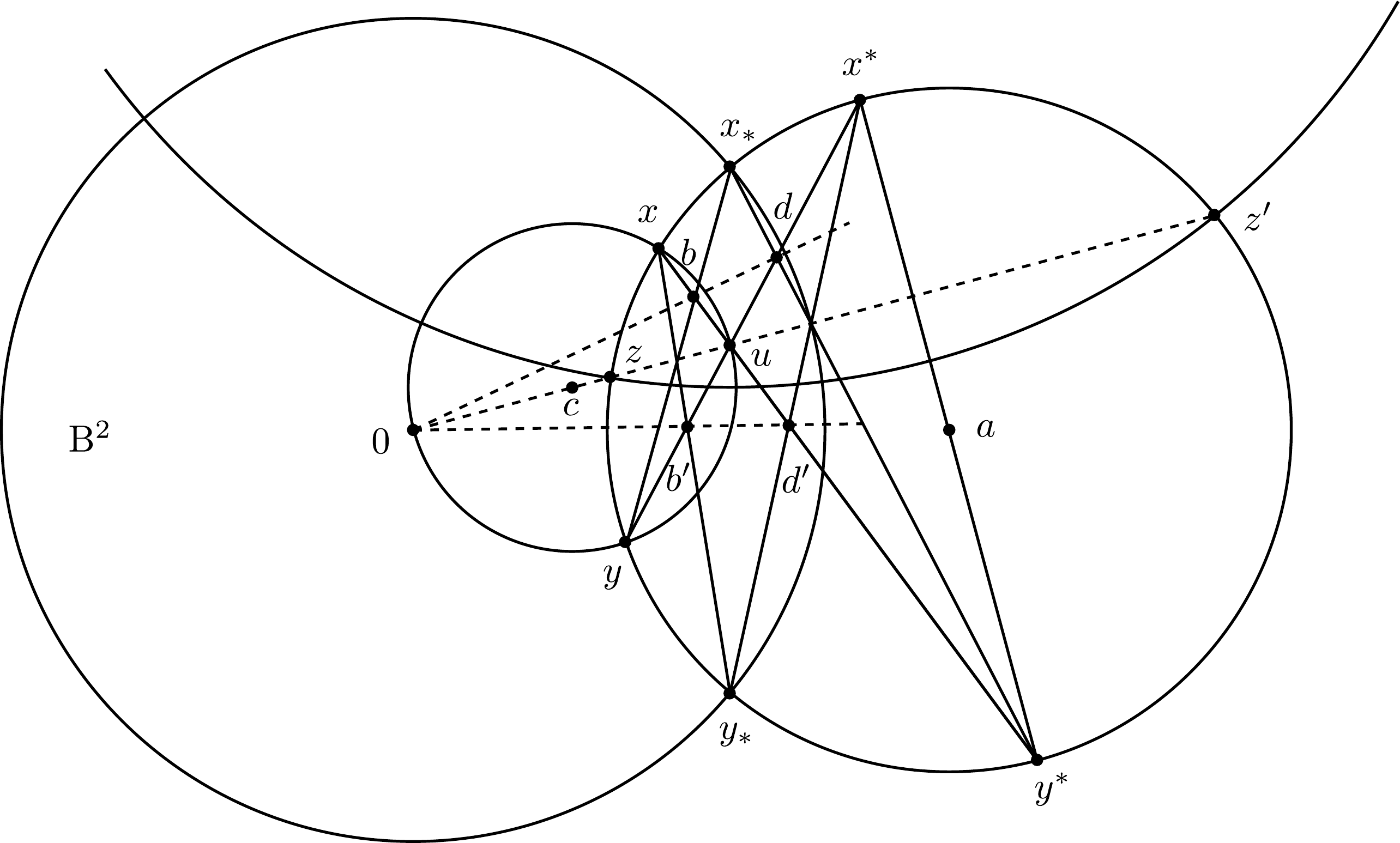}
\caption{\label{cor}}
\end{figure}

\begin{proposition} Let $x,y\in \BB\setminus\{0\}$ such that $0,x,y$ are noncollinear and $|x|\neq|y|$.
Let $x_*, y_*$, $u$, $ S^1(a,r_a), S^1(w,r_w)$ be as in Lemma 4.6,
$\sigma$ be the circular arc $x^*, x, y, y^*$. Construct the circle $S^1(c,r_c)$, which is through the points $0,x,y$, and let $z=L(0,c)\cap S^1(a,r_a)\cap \BB$, $z'=L(0,c)\cap S^1(a,r_a)\cap(\RR\setminus\BB)$. Let $b=L(x_*,y)\cap L(x,y^*)$, $d=L(x^*,y)\cap L(x_*,y^*)$, $b'=L(x,y_*)\cap L(x^*,y)$, $d'=L(x,y^*)\cap L(x^*,y_*)$, see Figure \ref{cor}.  Then

(1) $z$ is the midpoint of the hyperbolic segment $J[x,y]$, $S^1(c,r_c)$ is orthogonal to  $S^1(a,r_a)$ and $|z,x,x^*,z'|=|z,y,y^*,z'|$, provided $\sigma$ is a semicircle;

(2) the points $0, b, d$ are collinear, so are the points $0,b', d'$.
\end{proposition}

\begin{proof}(1) It is easy to see that $u$ is the orthocenter of the triangle $\triangle 0x^*y^*$ by the assumption of $\sigma$ and hence the segment $[0,u]$ is the diameter of $S^1(c,r_c)$. Therefore, $z$ is the midpoint of the hyperbolic segment $J[x,y]$ by Lemma 4.6(2).

By geometric observation, we obtain that
$$\angle cx0=\angle x0u=\angle xyx^*=\angle xy^*a=\angle axy^*,$$
therefore, $\angle cxa=\pi/2$ and $S^1(c,r_c)$ is orthogonal to  $S^1(a,r_a)$.

By Ptolemy's Theorem,
\begin{eqnarray}\label{zz'}
& &|z,x,x^*,z'|=|z,y,y^*,z'|\nonumber\\
&\Leftrightarrow& \frac{|z-x^*||x-z'|}{|z-x||x^*-z'|}=\frac{|z-y^*||y-z'|}{|z-y||y^*-z'|}\nonumber\\
&\Leftrightarrow& \frac{|x-x^*|}{|y-y^*|}=\frac{|z-x|}{|z'-y^*|}\cdot\frac{|z'-x^*|}{|z-y|}.
\end{eqnarray}
By similar triangles, we have $\frac{|x-x^*|}{|y-y^*|}=\frac{|x^*-u|}{|y^*-u|}$, $\frac{|z-x|}{|z'-y^*|}=\frac{|z-u|}{|y^*-u|}$, $\frac{|z'-x^*|}{|z-y|}=\frac{|x^*-u|}{|z-u|}$  and thus (\ref{zz'}) holds.

(2) By symmetry, we only need to prove that $0, b, d$ are collinear.

Without loss of generality, we may assume that $a_2=0$ and $x_*=w_1+i m$, where $w_1$ and $m$ are as in (\ref{w1w2a2=0}) and (\ref{ma1}), respectively.
First, we observe that $b=b_1+ib_2$ satisfies the following equations:
\beq\label{b1b2}
\left\{\begin{array}{ll}
b_2-x_2=\frac{y_2-x_2|y|^2}{y_1-x_1|y|^2}(b_1-x_1)\\
b_2-y_2=\frac{y_2-m}{y_1-w_1}(b_1-y_1)
\end{array}\right..
\eeq
Solving (\ref{b1b2}), we have
\beq\label{bb}
\frac{b_2}{b_1}=\frac{(y_1m-y_2w_1)(y_2-x_2|y|^2)+(x_1y_2-x_2y_1)(y_2-m)}{(y_1m-y_2w_1)(y_1-x_1|y|^2)+(x_1y_2-x_2y_1)(y_1-w_1)}.
\eeq
Similarly, the point $d=d_1+id_2$ satisfies the following equality
\beq\label{dd}
\frac{d_2}{d_1}=\frac{(y_1m-y_2w_1)(x_2-y_2|x|^2)-(x_1y_2-x_2y_1)(y_2-m|y|^2)}{(y_1m-y_2w_1)(x_1-y_1|x|^2)-(x_1y_2-x_2y_1)(y_1-w_1|y|^2)}.
\eeq

By (\ref{w1w2a2=0}), (\ref{y1-w1|y|^2}) and the assumption of $a_2=0$, we have that $0, b, d$ are collinear
\begin{eqnarray*}
&\Leftrightarrow& \frac{b_2}{b_1}=\frac{d_2}{d_1}\\
&\Leftrightarrow& (y_1m-y_2w_1)(y_2-x_2|y|^2)+(x_1y_2-x_2y_1)(y_2-m)\\
& &=(y_1m-y_2w_1)(x_2-y_2|x|^2)-(x_1y_2-x_2y_1)(y_2-m|y|^2)\\
&\Leftrightarrow& w_1(y_2-x_2+y_2|x|^2-x_2|y|^2)-m[y_1(1+|x|^2)-x_1(1+|y|^2)]=2(x_1y_2-x_2y_1)\\
&\Leftrightarrow& (y_2-x_2)[y_1(1+|x|^2)-x_1(1+|y|^2)]=0,
\end{eqnarray*}
this completes the proof.
\end{proof}

By the orthogonality of two circles, we easily obtain the following proposition.

\begin{proposition} Let $x,y\in \BB\setminus\{0\}$ such that $0,x,y$ are noncollinear. Let $0$ be the inversion point of points $x$ and $y$ with respect to the circles $S^1(x^*,t_x)$ and $S^1(y^*,t_y)$, respectively, where
$t_x=\sqrt{\frac{1}{|x|^2}-1}$ and $t_y=\sqrt{\frac{1}{|y|^2}-1}$. Then $S^1(x^*,t_x)$ is orthogonal to $S^1(y^*,t_y)$ if and only if $\cos\angle{x0y}=|x||y|$.
\end{proposition}

\begin{remark}It is easy to see that $x$ and $y$ are inversion points with respect to the circle $S^1(w,r_w)$, so are the pairs of $x_*, y_*$ and $x^*, y^*$. There are similar conclusions for other pairs of points with respect to the corresponding circles.
\end{remark}

\medskip
\noindent{\bf 4.36 Bisection of geodesic segment in $\BB$}
\medskip

We now provide several constructions of the midpoint of the hyperbolic geodesic segment in $\BB$. For this purpose, we only need to deal with two cases: the points $0, x, y$ are collinear and  $x,y$ are on a circle, which is orthogonal to $S^1$.

{\bf Case 1.} {\it The hyperbolic segment $J[x,y]$ is on the diameter of $\BB$.}

Without loss of generality we may assume that $\mbox{Im}x=\mbox{Im}y=0$.

\vspace{-3.2mm}
\begin{table}[htbp]
\centering
\begin{tabularx}{\textwidth}{lX}
\emph{Step(1)} & Construct the line $L(x,y)$.\\
\emph{Step(2)} & Construct the lines $L(x)$ and $L(y)$. Let $\{m,\overline{m}\}=L(x)\cap S^1$ and $\{n,\overline{n}\}=L(y)\cap S^1$.\\
\emph{Step(3)} & Construct the lines $L(m,\overline{n})$ and $L(\overline{m},n)$.
\end{tabularx}
\end{table}
\vspace{-3.8mm}

\begin{figure}[h]
\begin{minipage}[t]{0.45\linewidth}
\centering
\includegraphics[width=7cm]{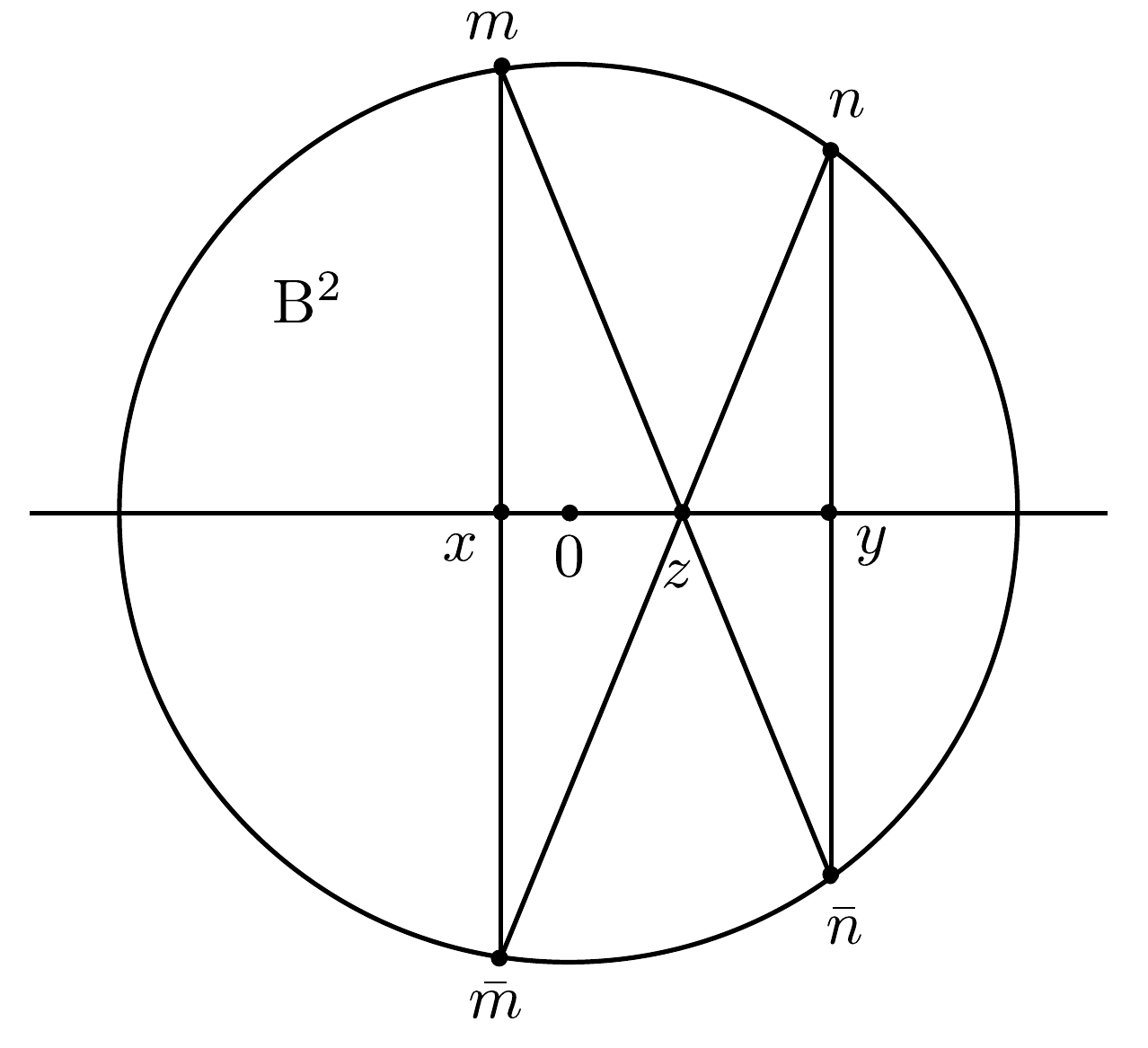}
\caption{\label{41}}
\end{minipage}
\hfill
\begin{minipage}[t]{0.45\linewidth}
\centering
\includegraphics[width=6.5cm]{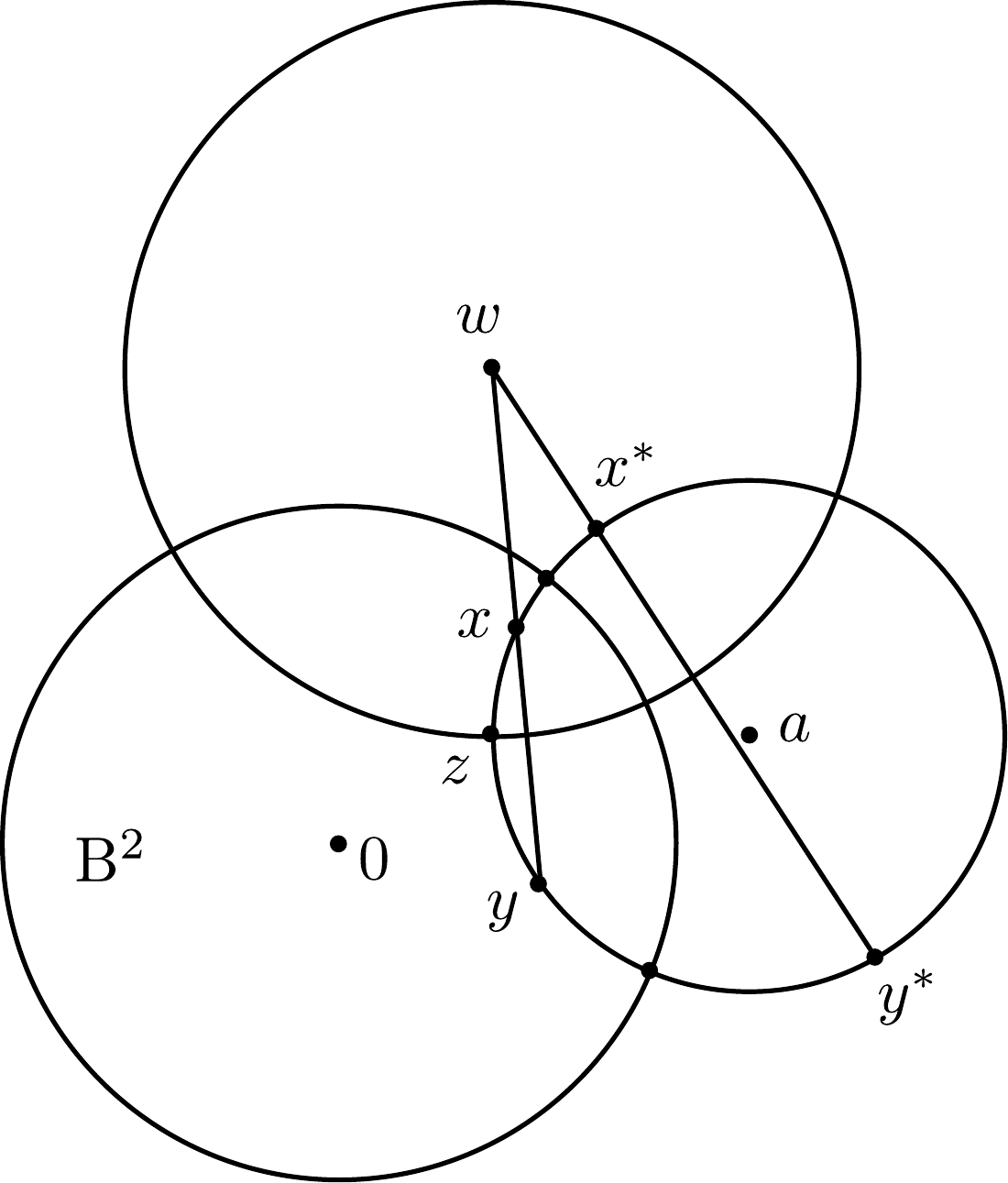}
\caption{\label{42}}
\end{minipage}
\end{figure}

\noindent Then the midpoint $z$ of $J[x,y]$ is the intersection of $L(m,\overline{n})$ and $L(\overline{m},n)$ by Lemma 3.1(4) and Projection property 4.1, see Figure \ref{41}.

{\bf Case 2.} {\it The hyperbolic segment $J[x,y]$ is not on the diameter of $\BB$. That means $J[x,y]$ is on a circle that is orthogonal to $S^1$.}

If $|x|=|y|$, then by (\ref{rhowv}) it is clear that the midpoint of $J[x,y]$ is the intersection of $L(0,a)$ and $S^1(a,r_a)$ in $\BB$ which is orthogonal to $S^1$. Therefore we assume that $|x|\neq|y|$ in the sequel.

{\noindent\bf Method I.}\\
\emph{Step (1)} Construct the circle $S^1(a,r_a)$ which contains the four points $x,y,x^*,y^*$.\\
\emph{Step (2)} Construct the lines $L(x,y)$ and $L(x^*,y^*)$. Let
$w=L(x,y)\cap L(x^*,y^*)$.\\
\emph{Step (3)} Construct the circle  $S^1(w,r_w)$ which is orthogonal to the circle $S^1(a,r_a)$.\\
Then the midpoint $z$ of $J[x,y]$ is the intersection of $S^1(w,r_w)$ and $S^1(a,r_a)$ in $\BB$ by Lemma 4.6(1), see Figure \ref{42}.

\medskip
\noindent{\bf Remark 4.37} {\it Method I is similar  to the construction in \cite[2.9]{kv}.}
\medskip

{\noindent\bf Method II.-VI.}

\vspace{-3.2mm}
\begin{table}[htbp]
\centering
\begin{tabularx}{\textwidth}{lX}
\emph{Step (1)} & Construct the circle $S^1(a,r_a)$ which contains the four points $x,y,x^*,y^*$.
Let $\{x_*,y_*\}=S^1(a,r_a)\cap S^1$, $x_*,x,y,y_*$ occur in this order on $S^1(a,r_a)$.\\
\emph{Step (2)}
\end{tabularx}
\end{table}
\vspace{-4mm}

II. Construct the lines $L(x,y^*)$ and $L(y,x^*)$, let $u=L(x,y^*)\cap L(y,x^*)$.

III. Construct the lines $L(x,x_*)$ and $L(y,y_*)$, let $v=L(x,x_*)\cap L(y,y_*)$.

IV.  Construct the lines $L(x,y_*)$ and $L(y,x_*)$, let $s=L(x,y_*)\cap L(y,x_*)$.

V. Construct the lines $L(x_*,y^*)$ and $L(y_*,x^*)$, let $t=L(x_*,y^*)\cap L(y_*,x^*)$.

\begin{figure}[h]
\begin{minipage}[t]{0.45\linewidth}
\centering
\includegraphics[width=6.5cm]{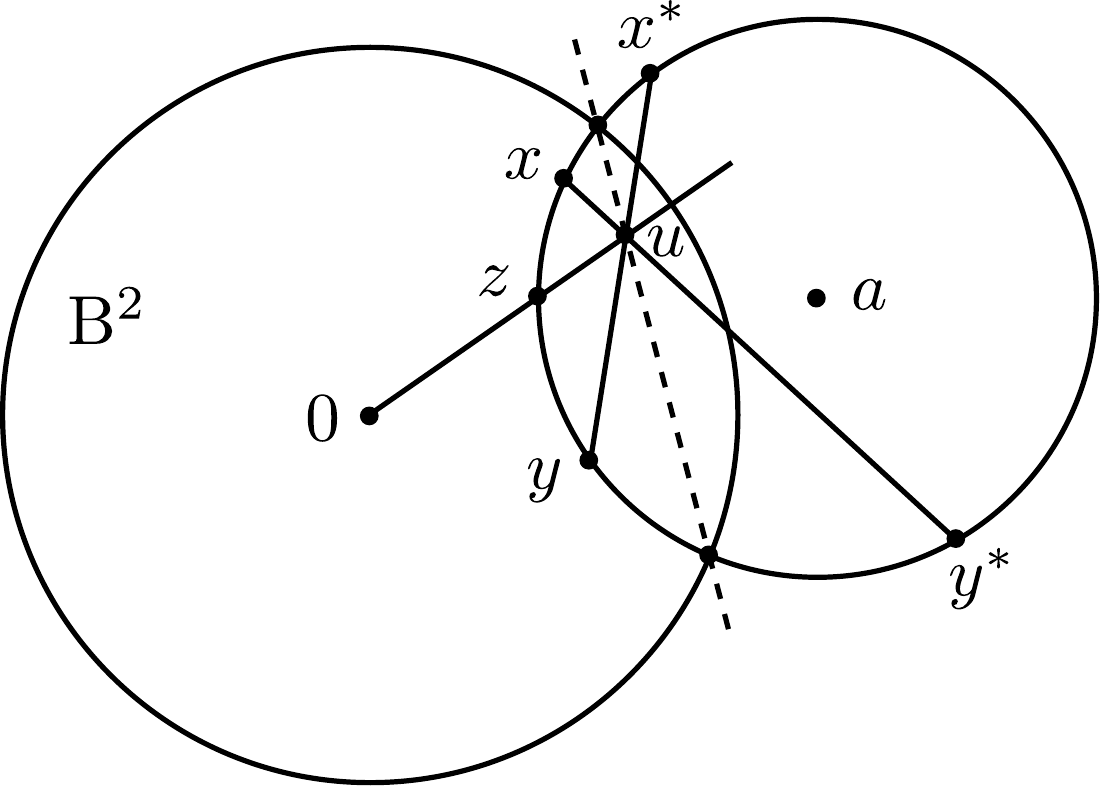}
\caption{\label{43}}
\end{minipage}
\hfill
\begin{minipage}[t]{0.45\linewidth}
\centering
\includegraphics[width=6.5cm]{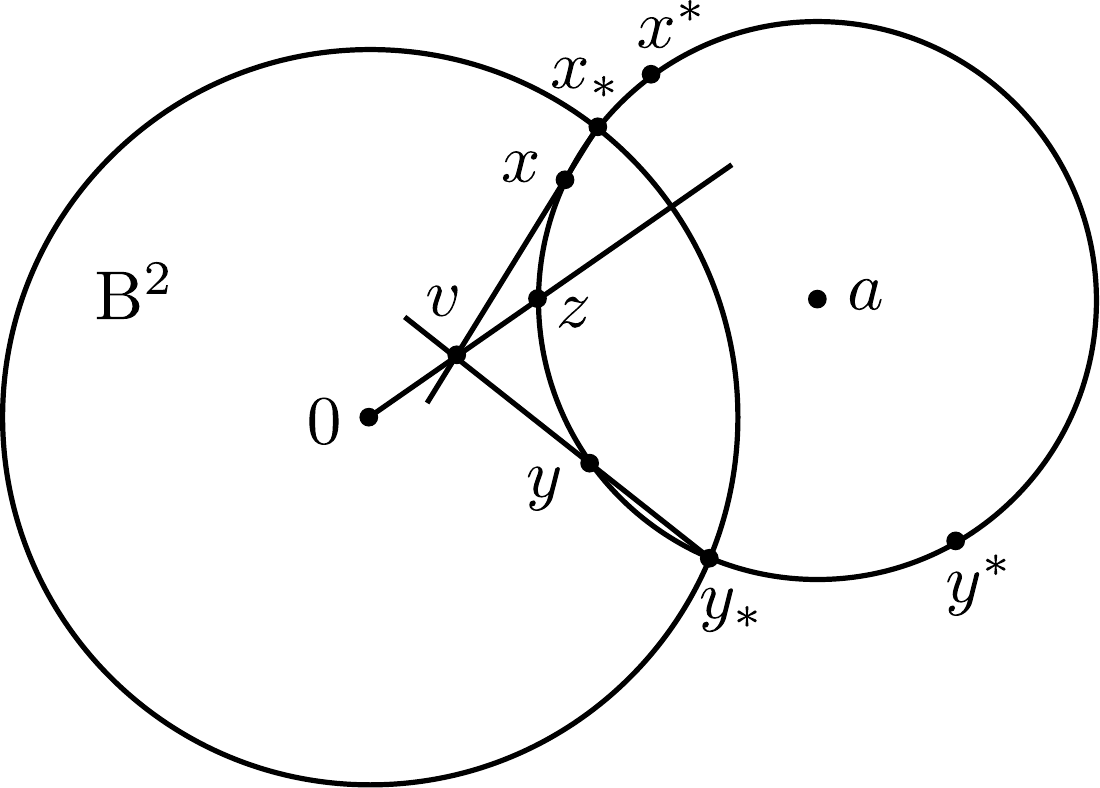}
\caption{\label{44}}
\end{minipage}
\end{figure}

\begin{figure}[h]
\begin{minipage}[t]{0.45\linewidth}
\centering
\includegraphics[width=6.5cm]{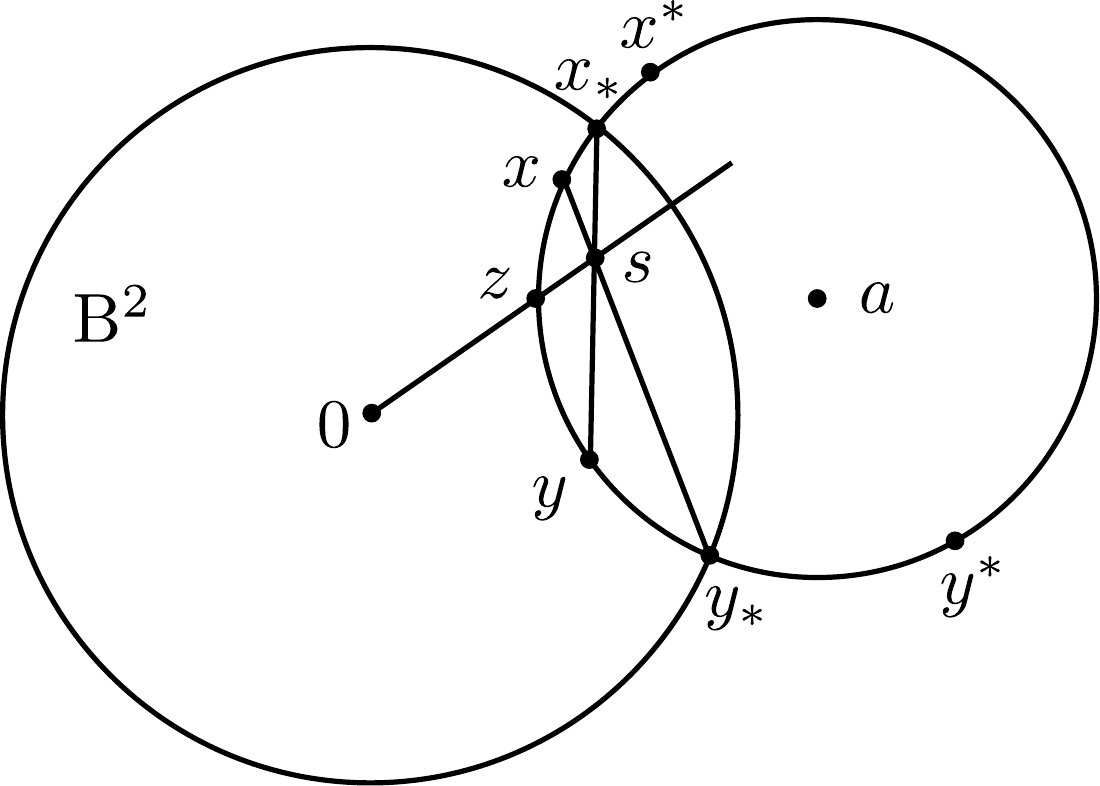}
\caption{\label{45}}
\end{minipage}
\hfill
\begin{minipage}[t]{0.45\linewidth}
\centering
\includegraphics[width=6.5cm]{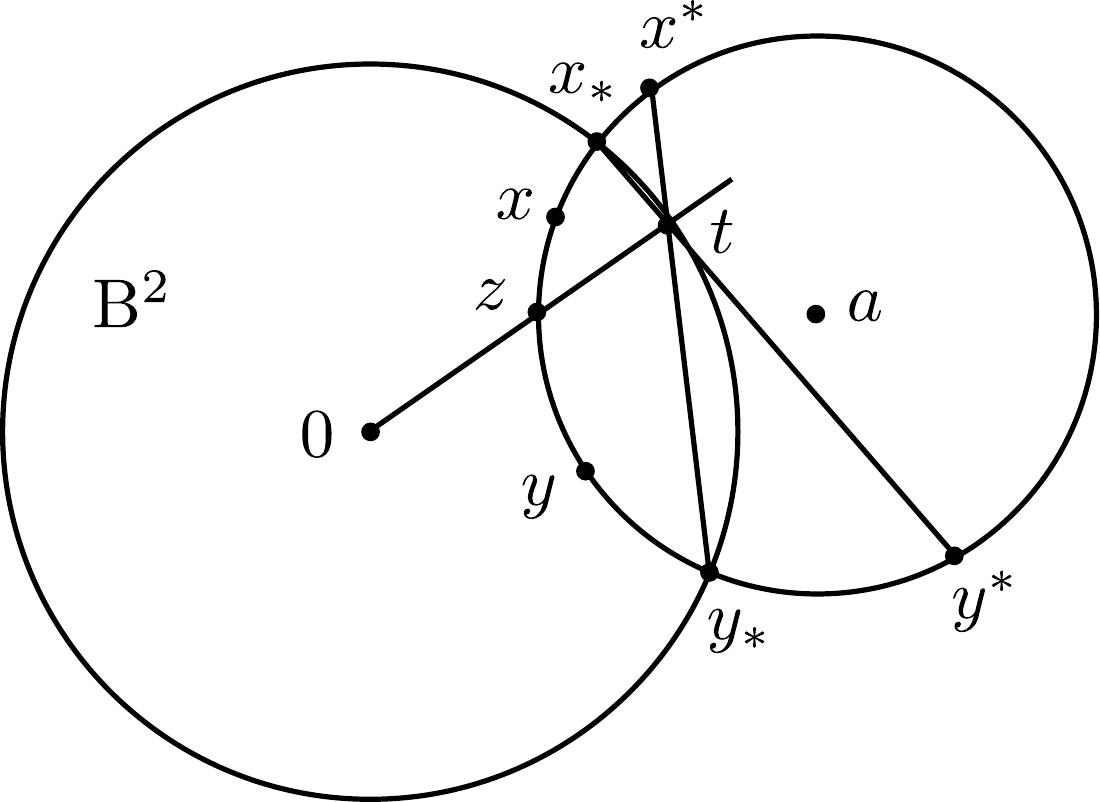}
\caption{\label{46}}
\end{minipage}
\end{figure}

\begin{figure}[h]
\begin{minipage}[t]{0.6\linewidth}
\centering
\includegraphics[width=8cm]{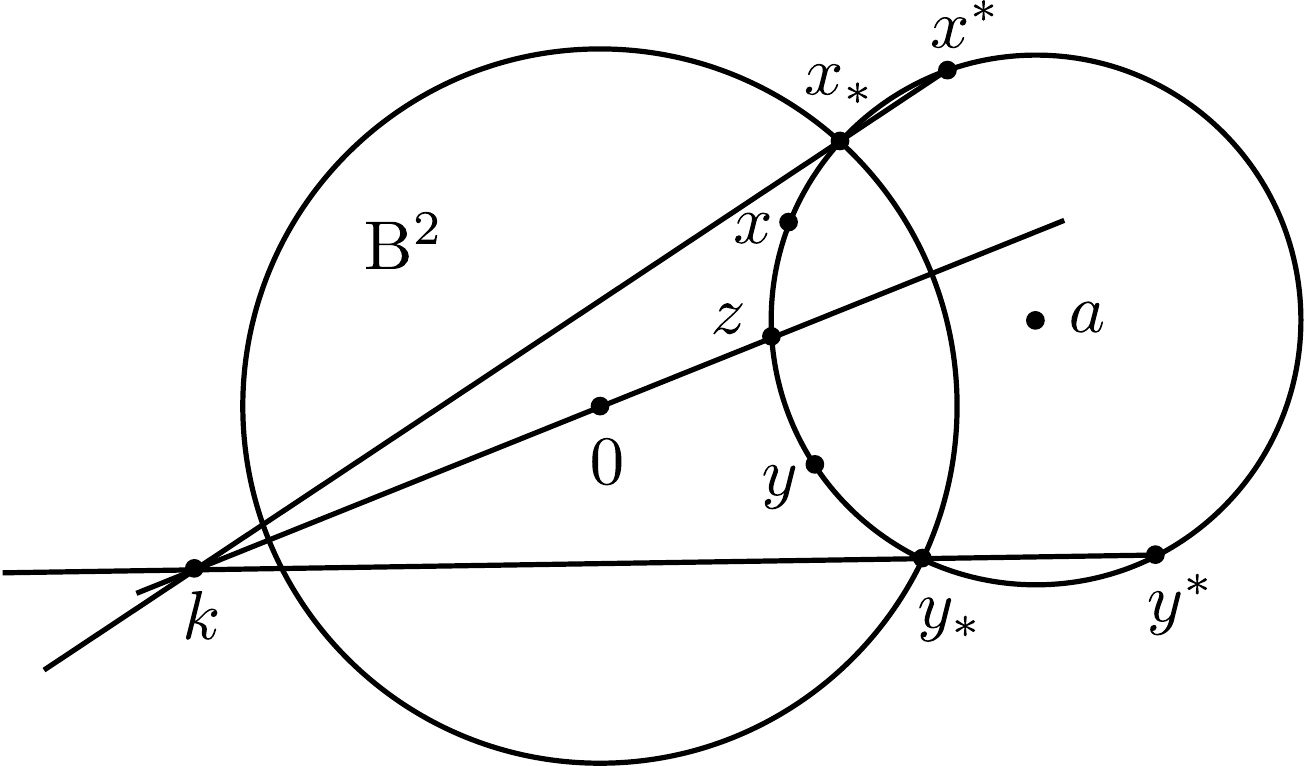}
\caption{\label{47}}
\end{minipage}
\end{figure}

VI. Construct the lines $L(x_*,x^*)$ and $L(y_*,y^*)$, let $k=L(x_*,x^*)\cap L(y_*,y^*)$.\\
\emph{\,\,\,Step (3)} Construct the line $L(0,g)$ for each $g\in\{u,v,s,t,k\}$.\\
Then the midpoint $z$ of $J[x,y]$ is the intersection  of $L(0,g)$ and $S^1(a,r_a)$ in $\BB$ by Lemma 4.6(2), see Figures \ref{43}--\ref{47}.

\subsection*{Acknowledgments}
The research of Matti Vuorinen was supported by the Academy of Finland,
Project 2600066611. The research of Gendi Wang was supported by CIMO
of Finland, Grant TM-10-7364. The authors thank the referee for a valuable set of corrections.

\end{document}